\documentclass[11pt]{article} 
 
\usepackage{latexsym} 
\usepackage{amsfonts} 
\usepackage{amsmath} 
\usepackage{amsthm} 
\usepackage{mathrsfs} 
\usepackage{dsfont} 
\usepackage{bbold} 
\usepackage[english]{babel} 
\usepackage{caption} 
\usepackage{epsfig} 
\usepackage{mathtools} 
\usepackage{enumitem}
\usepackage{float} 
\usepackage{subfigure} 
\usepackage{psfrag} 
\usepackage{graphicx} 
\usepackage{epsfig} 
\usepackage{amsbsy} 
\usepackage{hyperref}
 
\DeclareMathOperator{\Val}{\matV}

\newtheorem{theorem}{Theorem} 
\newtheorem*{prop*}{Theorem}

\newtheorem{defi}[theorem]{Definition} 
\newtheorem{lemma}[theorem]{Lemma} 
\newtheorem{prop}[theorem]{Proposition} 
\newtheorem{rmk}[theorem]{Remark}

\newcommand{\zerarcounters}{\setcounter{equation}{0}\setcounter{theorem}{0}}

\newcommand{\ZZZ}{\mathds{Z}} 
\newcommand{\CCC}{\mathds{C}} 
\newcommand{\NNN}{\mathds{N}} 
 
\newcommand{\RRR}{\mathds{R}} 
\newcommand{\TTT}{\mathds{T}}

\newcommand{\calA}{{\mathcal A}} 
\newcommand{\calB}{{\mathcal B}} 
\newcommand{\calC}{{\mathcal C}} 
\newcommand{\calD}{{\mathcal D}} 
 
\newcommand{\calF}{{\mathcal F}} 
\newcommand{\calG}{{\mathcal G}}

\newcommand{\LL}{{\mathcal L}}

\newcommand{\calP}{{\mathcal P}} 
 
\newcommand{\RR}{{\mathcal R}} 
 
\newcommand{\TT}{{\mathcal T}} 
\newcommand{\calU}{{\mathcal U}}

\newcommand{\gotB}{{\mathfrak B}} 
\newcommand{\gotC}{{\mathfrak C}}

\newcommand{\gotF}{{\mathfrak F}} 
\newcommand{\gotG}{{\mathfrak G}}

\newcommand{\gotS}{{\mathfrak S}} 
\newcommand{\gotT}{{\mathfrak T}}

\newcommand{\matV}{{\mathscr V}}

\newcommand{\prova}{\noindent{\it Proof. }} 
\newcommand{\io}{\infty} 
\newcommand{\e}{\varepsilon} 
\newcommand{\al}{\alpha} 
\newcommand{\de}{\delta} 
\newcommand{\be}{\beta} 
\newcommand{\n}{\nu} 
 
\newcommand{\x}{\xi}

\newcommand{\om}{\omega}

\newcommand{\s}{\sigma} 
\newcommand{\ols}{\overline{\sigma}}

\newcommand{\oo}{\boldsymbol{\omega}}

\newcommand{\nn}{\boldsymbol{\nu}}

\newcommand{\bt}{\boldsymbol{t}}

\newcommand{\ii}{{\rm i}}

\def\ins#1#2#3{\vbox to0pt{\kern-#2 \hbox{\kern#1 #3}\vss}\nointerlineskip} 
 
\begin{document}
 
%%%%%%%%%%%%%%%%%%%%%%%%%%%%%%%%%%%%%%%%%%%%%%%%%
\title{\bf Invariant curves for exact symplectic twist maps\\
of the cylinder with Bryuno rotation numbers}
%%%%%%%%%%%%%%%%%%%%%%%%%%%%%%%%%%%%%%%%%%%%%%%%%
 
\author {\bf Guido Gentile\vspace{2mm} \\ 
\small Dipartimento di Matematica e Fisica, Universit\`a Roma Tre, Roma, I-00146, Italy \\ 
\small Email: gentile@mat.uniroma3.it}
\date{} 
 
\maketitle 
 
%%%%%%%%%%%%%%%%%%%%%%%%%%%%%%%%%%%%%%%%%%%%%%%%%
\begin{abstract} 
Since Moser's seminal work it is well known that the invariant curves of smooth
nearly integrable twist maps of the cylinder with Diophantine rotation number 
are preserved under perturbation. In this paper we show that, in the analytic class, 
the result extends to Bryuno rotation numbers. First, we will show that
the series expansion for the invariant curves in powers of 
the perturbation parameter can be formally defined, then we shall prove
that the series converges absolutely in a neighbourhood of the origin.
This will be achieved using multiscale analysis and renormalisation group techniques
to express the coefficients of the series as sums of values which are represented
graphically as tree diagrams and then exploit cancellations
between terms contributing to the same perturbation order.
As a byproduct we shall see that, when perturbing linear maps, the series expansion 
for an analytic invariant curve converges  for all perturbations 
if and only if the corresponding rotation number satisfies the Bryuno condition.
\end{abstract} 
%%%%%%%%%%%%%%%%%%%%%%%%%%%%%%%%%%%%%%%%%%%%%%%%%
  
%%%%%%%%%%%%%%%%%%%%%%%%%%%%%%%%%%%%%%%%%%%%%%%%%
%%%%%%%%%%%%%%%%%%%%%%%%%%%%%%%%%%%%%%%%%%%%%%%%%
\zerarcounters 
\section{Introduction} 
\label{sec:1} 
%%%%%%%%%%%%%%%%%%%%%%%%%%%%%%%%%%%%%%%%%%%%%%%%%
%%%%%%%%%%%%%%%%%%%%%%%%%%%%%%%%%%%%%%%%%%%%%%%%%

Consider an exact symplectic twist map $\Phi$ on the cylinder $\TTT\times\RRR$
\begin{equation} \label{eq:1.1}
\begin{cases}
x' = x + a(y) + \e \, f(x,y) , \\ y'=y + \e \, g(x,y) ,
\end{cases}
\end{equation}
where $\e\in\RRR$ and the functions $a,f,g$ are analytic,
with $a'(y)>0$ (\emph{twist condition}), and $2\pi$-periodic in $x$.
To fix notations we shall assume that the functions can be analytically
extended to a complex domain $\calD:=\TTT_{\xi_0}\times\calA$,
where $\TTT_{\xi_0}=\{x\in\CCC/2\pi\ZZZ : |\hbox{Im}\,x| \le \xi_0\}$ and
$\calA\subset\CCC$ is a complex neighbourhood of an interval of the real axis.
We are interested in the case where the supremum norms of $f$ and $g$ satisfy
$\max\{|f|,|g|\}=1$ in $\calD$ and $|\e| \le \e_0$, with $\e_0$ small enough.
For that reason, $\e$ is called the \emph{perturbation parameter}.

For $\e=0$ the map is integrable: $y$ is kept constant, $y=y_0$, and the dynamics is
a rotation of the variable $x$ by an angle $2\pi\om=a(y_0)$.
Below, following Moser's example \cite{M1}, we will consider initially the system
\begin{equation} \label{eq:1.2}
\begin{cases}
x' = x + y + \e \, f(x,y) , \\ y'=y + \e \, g(x,y) 
\end{cases}
\end{equation}
and eventually discuss how the analysis has to be changed to deal with
the general case (\ref{eq:1.1}).

Since we are assuming that $\Phi$ is exact symplectic, there exists a generating function
$S_{\e}(x,x')=S_0(x,x')+\e \s_{\e}(x,x')$ such that (\ref{eq:1.1}) can be written in the form
\begin{equation} \label{eq:1.3}
y = - \frac{\partial S_{\e}(x,x')}{\partial x} , \qquad
y' = \frac{\partial S_{\e}(x,x')}{\partial x'} .
\end{equation}
For the system (\ref{eq:1.2}) one has
\begin{equation} \nonumber
S_{0}(x,x')= \frac{1}{2} (x'-x)^{2} 
\end{equation}
In that case, if we take $\sigma_{\e}(x,x')=1-\cos x$,
so that $f(x,y)=g(x,y)=\sin x$ in (\ref{eq:1.2}),
we recover the Taylor-Chirikov \emph{standard map}, a map widely studied in the literature
since the original papers by Chirikov and Greene \cite{Ch,Gr,LL}.

Requiring that the map $\Phi$ is exact symplectic is equivalent to require that
it is area-preserving and has zero flux: in particular the latter condition
is automatically satisfied if $\Phi$ is area preserving and preserves the boundaries
of an annulus \cite{Go}. In the following it will be more convenient to study
the system (\ref{eq:1.1}) in the form (\ref{eq:1.3})
to better exploit the symplectic nature of the map.

We say that the map (\ref{eq:1.1}) admits an invariant curve if there exist two
functions $h,H$ of period $2\pi$ such that $x=\psi+h(\psi)$, $y=y_{0}+H(\psi)$ and the
dynamics induced on the curve is given by $\psi'=\psi+2\pi\om$, with $2\pi\om=a(y_{0})$.
The number $\om$ will be called the \emph{rotation number} of the curve.
If the functions $h$ and $H$ are analytic, we say that the invariant curve is analytic.
The existence of invariant curves for exact symplectic twist maps on the cylinder
was first proved by Moser in the differentiable class
by assuming a standard Diophantine condition on $\om$ \cite{M1}.
Moreover the result extends to more general maps, satisfying
the intersection property \cite{SM}: the latter means that every smooth curve
wrapping around the cylinder intersects itself under the action of the map.
However, in the analytic class the Diophantine condition is not expected to be optimal.
In the present paper we show that one can be take a weaker Diophantine condition
on the rotation number $\om$: the \emph{Bryuno condition} \cite{B},
defined as $\gotB(\om)<\io$, where
\begin{equation} \label{eq:1.4}
\gotB(\om) := \sum_{k=0}^{\io} \frac{1}{q_{k}} \log q_{k+1} 
\end{equation}
and $\{q_{k}\}$ are the denominators of the continued fraction expansion of $\om$ \cite{S}. 
In fact, we prove the following result.

%%%%%%%%%%%%%%%%%%%%%%%%%%%%%%%%%%%%%%%%%%%%%%%%%
\begin{theorem} \label{thm:1.1}
Consider an exact symplectic map on the cylinder of the form (\ref{eq:1.1}),
with $f$ and $g$ analytic. If $\omega$ satisfies the Bryuno condition
then there exists $\e_{0}>0$ such that the map admits
an analytic invariant curve with rotation number $\omega$
analytic in $\e$ for $|\e|<\e_{0}$.
\end{theorem}
%%%%%%%%%%%%%%%%%%%%%%%%%%%%%%%%%%%%%%%%%%%%%%%%%

This will be proved in Section \ref{sec:2} to \ref{sec:5} for maps of the form (\ref{eq:1.2})
and then extended to any maps of the form (\ref{eq:1.1}) in Section \ref{sec:6}.
As a byproduct of the proof we show directly that the perturbation series
for the invariant curve converges absolutely. This requires some work,
since the standard Lindstedt algorithm produces the coefficients as sums
of many terms which can grow proportionally to factorials, so one has to
show that there are compensations between such terms. We shall explicitly
exhibit the cancellations which make the sum of all terms contributing to
the same perturbation order $k$ bounded by some constant to the power $k$.
The cancellation mechanisms is similar to that found in the case
of Hamiltonian flows \cite{E,Ga,GM}, even though the implementation
is sligthly different. Unlike the case of flows we shall be able to obtain a closed
equation involving only the conjugation function of the angle variable
(which is not the standard approach for flows, even though there are remarkable
exceptions \cite{M3,SZ}). On the other hand this will make the cancellation analysis
a bit more delicate. Moreover, in the general case (\ref{eq:1.1}), the cancellations
will be showed to involve identities that can be of interest in their own;
see in particular Appendix \ref{app:b}.

Recently the Bryuno condition received a lot of attention in the literature
on small divisor problems arising in the context of analytic dynamical systems, including
circle diffeomorphisms \cite{Y2}, holomorphic maps on the plane \cite{Y1,BC,GiM},
maximal and lower-dimensional tori in quasi-integrable systems \cite{Ge3,KK1,KK2},
skew-product systems \cite{Ge2,CM}. The Bryuno condition for twist maps of the cylinder
has been explicitly considered in the literature only in special cases, such as \cite{BG}.
Indeed, most of the results on the existence of invariant tori under the Bryuno condition
are for vector fields. In general a quasi-integrable discrete map can be interpolated
by a continuous flow with the same number of degrees of freedom only approximately,
even though the error of the approximation can be made exponentially small
in the perturbation parameter \cite{De,N,BeGi}. However,
the existence of invariant curves be deduced directly from the results for continuous systems
by considering non-autonomous Hamiltonian systems whose flows interpolate
the maps at discrete times \cite{KP} (see also \cite{M2,SZ} in the case of smooth,
not necessarily nearly integrable maps). Of course, this a way to derive the results.
Actually, the construction of the interpolating flow can be intricate:
for instance, a simple-looking system such as the standard map naturally leads
to a singular Hamiltonian and already a non-analytic smoothing yields a much less handy system \cite{M2}.
Therefore we think that it could be interesting to have a self-contained proof in the case of twist maps.

The Bryuno condition is known to be optimal in the case of the Siegel problem \cite{Y2} and
the circle diffeomorphisms \cite{Y1}. It has been conjectured to be optimal also in the case
of twist maps by MacKay \cite{Mac}. However, despite being supported by numerical investigations \cite{MS},
the conjecture remains still unproved. The main result in this direction has been obtained by Forni \cite{F}
(see also \cite{Be} for extensions to higher-dimensional flows), who proved that
if $\Phi_{0}$ is an integrable twist map and $\om$ satisfies a condition
somewhat stronger than the \emph{converse Bryuno condition} $\gotB(\om)=\io$, more precisely
\begin{equation} \label{eq:1.5}
\limsup_{n\to +\io} \frac{1}{q_{n}} \log q_{n+1} > 0 ,
\end{equation}
then there are perturbations arbitrarily close to $\Phi_{0}$ (in the analytic topology) which do 
not admit any invariant curve with rotation number $\omega$. As noticed in \cite{F},
the result does not immediately applies to one-parameter families of maps of the
form \eqref{eq:1.1} with the functions $f$ and $g$ either fixed or analytic in a fixed domain.
Here we show that, in the case of exact symplectic twist maps
of the cylinder of the form (\ref{eq:1.2}), i.e. perturbations of the linear map,
the Bryuno condition is optimal in order to have an invariant
curve depending analytically on $\e$ for any $f$ and $g$ analytic in
a domain $\calD$ as after \eqref{eq:1.1}. As a byproduct
analyticity in $\e$ automatically implies analyticity of the functions $h,H$ in $\psi$ and hence
analyticity of the invariant curve. More precisely we obtain the following result.

%%%%%%%%%%%%%%%%%%%%%%%%%%%%%%%%%%%%%%%%%%%%%%%%%
\begin{theorem} \label{thm:1.2}
Given the integrable map $\Phi_{0}$ on $\TTT\times\RRR$
\begin{equation} \nonumber
\begin{cases}
x' = x + y , \\ y'= y ,
\end{cases}
\end{equation}
let $\Phi$ be an analytic exact symplectic map of the form $\Phi=\Phi_{0}+\e\Phi_{1}$, 
with $\Phi_{1}=(f,g)$ analytic in a domain $\calD$ and such that
$|\Phi_{1}|:=\max\{|f|,|g|\}=1$ in $\calD$. Then $\om\in\RRR$ satisfies the Bryuno condition
if and only if there exists $\e_0>0$ such that any $\Phi$ admits an analytic invariant curve
with rotation number $\om$ analytic in $\e$ for $|\e|<\e_0$
\end{theorem}
%%%%%%%%%%%%%%%%%%%%%%%%%%%%%%%%%%%%%%%%%%%%%%%%%

The proof, to be given in Section \ref{sec:5}, which we refer to for more details, follows 
immediately from Theorem \ref{thm:1.1} and some results 
available in the literature \cite{D,BG}.

%%%%%%%%%%%%%%%%%%%%%%%%%%%%%%%%%%%%%%%%%%%%%%%%%
%%%%%%%%%%%%%%%%%%%%%%%%%%%%%%%%%%%%%%%%%%%%%%%%%
\zerarcounters 
\section{Formal power series and tree expansion}
\label{sec:2} 
%%%%%%%%%%%%%%%%%%%%%%%%%%%%%%%%%%%%%%%%%%%%%%%%%
%%%%%%%%%%%%%%%%%%%%%%%%%%%%%%%%%%%%%%%%%%%%%%%%%

We start by considering exact symplectic maps of the form (\ref{eq:1.2})
and defer to Section \ref{sec:6} how to extend the analysis to (\ref{eq:1.1}).
So, we prove first the following result, which is a special case of Theorem \ref{thm:1.1}.

%%%%%%%%%%%%%%%%%%%%%%%%%%%%%%%%%%%%%%%%%%%%%%%%%
\begin{theorem} \label{thm:2.0}
Consider an exact symplectic map on the cylinder of the form (\ref{eq:1.2}),
with $f$ and $g$ analytic. If $\omega$ satisfies the Bryuno condition
then there exists $\e_{0}>0$ such that the map admits
an analytic invariant curve with rotation number $\omega$
analytic in $\e$ for $|\e|<\e_{0}$.
\end{theorem}
%%%%%%%%%%%%%%%%%%%%%%%%%%%%%%%%%%%%%%%%%%%%%%%%%

As said in Section \ref{sec:1} we rewrite (\ref{eq:1.2})  in the form
\begin{equation} \label{eq:2.1}
\begin{dcases}
y' = x'-x + \e \, \frac{\partial \sigma_{\e}(x,x')}{\partial x'} , \\
y = x'-x - \e \, \frac{\partial \sigma_{\e}(x,x')}{\partial x} ,
\end{dcases}
\end{equation}
where $\sigma_{\e}(x,x')$ is a function analytic in $x$, $x'$ and $\e$,
such that $\sigma_{\e}(x+2\pi,x'+2\pi)=\sigma_{\e}(x,x')$.
It is more convenient to write $\sigma_{\e}(x,x')=\ols_{\e}(x,x-x')$,
with $\ols_{\e}(x,z)=\ols_{\e}(x+2\pi\om,z)$. For simplicity, we first consider the case
in which $\ols_{\e}(x,z)=\ols(z,z)$ does not depend on $\e$ and show at the end that
the analysis easily extends so as to take into account also the dependence on $\e$.

By the assumption of analyticity one has
\begin{equation} \label{eq:2.2}
\ols(x,z)=\sum_{\nu\in\ZZZ} {\rm e}^{\ii\nu x} \, \s_{\nu}(z) , \qquad
\frac{1}{q!} \left| \partial_{z}^{q} \s_{\nu}(2\pi\om) \right| \le \Xi \, {\rm e}^{-\xi |\nu|}\rho^{-q} \quad
\forall (\nu,q) \in\ZZZ\times\ZZZ_{+} ,
\end{equation}
for suitable positive constants $\Xi$, $\xi$ and $\rho$, 
depending on the domain $\calD$ introduced after (\ref{eq:1.1}). Here and henceforth
$\partial_{z}^q\s_{\n}(2\pi\om)$ stands for $\partial_{z}^q\s_{\n}(z)|_{z=2\pi\om}$.

We look for two analytic functions $h$ and $H$ (\emph{conjugation}) such that
\begin{equation} \label{eq:2.3}
x = \psi + h(\psi) , \qquad y=2\pi\om+H(\psi) 
\end{equation}
and the dynamics in terms of $\psi$ becomes $\psi'=\psi+2\pi\om$.
By introducing (\ref{eq:2.3}) into (\ref{eq:2.1}) and imposing $\psi'=\psi+2\pi\om$
one finds the functional equation for $h$
\begin{eqnarray} \label{eq:2.4}
& & h(\psi+2\pi\om) + h(\psi-2\pi\om)-2h(\psi) \nonumber \\
& & \qquad =
\e \, \partial_1 \ols(\psi+h(\psi),2\pi\om+h(\psi+2\pi\om)-h(\psi)) \\
& & \qquad \qquad 
- \, \e \, \partial_2 \ols(\psi+h(\psi),2\pi\om+h(\psi+2\pi\om) -h(\psi)) \nonumber \\
& & \qquad \qquad
+ \, \e \, \partial_2 \ols(\psi-2\pi\om+h(\psi-2\pi\om),2\pi\om+h(\psi)-h(\psi-2\pi\om)) , \nonumber
\end{eqnarray}
where $\partial_j$ denotes derivative with respect to the $j$-th argument.
If the equation (\ref{eq:2.4}) can be solved then the function $H$ can be 
obtained from
\begin{equation} \label{eq:2.5}
H(\psi) = h(\psi) - h(\psi-2\pi\om) + \e \,
\partial_2 \ols(\psi-2\pi\om+h(\psi-2\pi\om),2\pi\om+h(\psi)-h(\psi-2\pi\om)) .
\end{equation}
Moreover $H$ inherits the same regularity properties of $h$,
so that we can confine ourselves to the functional equation (\ref{eq:2.4}).

If there exists a solution to (\ref{eq:2.4}) we say that the system
admits an \emph{invariant curve} with rotation number $\om$; 
if the functions are analytic in $\psi$ and $\e$ we say that the invariant curve
is analytic in $\psi$ and $\e$.
In the following we will be interested in invariant curves of this kind.

By taking the formal Fourier-Taylor expansion of $h$,
\begin{equation} \label{eq:2.6}
h(\psi) = \sum_{k=1}^{\io} \e^{k}
\sum_{\nu\in\ZZZ} {\rm e}^{i\nu\psi} h^{(k)}_{\nu} ,
\end{equation}
one can rewrite (\ref{eq:2.4}) as infinitely many equations
\begin{eqnarray} \label{eq:2.7}
& & \hskip-.5truecm \de(\om\nu) \, h^{(k)}_{\nu} =
\sum_{p,q\ge 0} 
\sum_{\substack{\nu_0+\nu_{1}+\ldots+\nu_{p}\\+\mu_{1}+\ldots+\mu_{q}=\nu}}
\sum_{\substack{k_{1}+\ldots+k_{p}\\+k'_{1}+\ldots+k'_{q}=k-1}}
\frac{1}{p!q!} \, (\ii\nu_0)^{p+1} \partial_{z}^{q} \s_{\nu_0}(2\pi\om) \times \nonumber \\
& & \quad  
\times \Big( \prod_{i=1}^{q} \de_{+}(\om\mu_{i}) \Big)
h^{(k_1)}_{\nu_1}\ldots h^{(k_p)}_{\nu_p}
h^{(k'_1)}_{\mu_1}\ldots h^{(k'_q)}_{\mu_q} \\
& & \qquad \quad + \sum_{p,q\ge 0} 
\sum_{\substack{\nu_0+\nu_{1}+\ldots+\nu_{p}\\+\mu_{1}+\ldots+\mu_{q}=\nu}}
\sum_{\substack{k_{1}+\ldots+k_{p}\\+k'_{1}+\ldots+k'_{q}=k-1}}
\frac{1}{p!q!} \, (\ii\nu_0)^{p} \partial_{z}^{q+1} \s_{\nu_0}(2\pi\om) \times \nonumber \\
& & \quad 
\times \left( {\rm e}^{-2\pi\ii \om(\nu_0+\nu_1+\ldots+\nu_p)}
\Big( \prod_{i=1}^{q} \de_{-}(\om\mu_{i}) \Big) -
\Big( \prod_{i=1}^{q} \de_{+}(\om\mu_{i}) \Big) \right)
h^{(k_1)}_{\nu_1}\ldots h^{(k_p)}_{\nu_p}
h^{(k'_1)}_{\mu_1}\ldots h^{(k'_q)}_{\mu_q} , \nonumber
\end{eqnarray}
where we have set 
\begin{equation} \label{eq:2.0}
\de(u):=2(\cos 2\pi u -1) , \qquad
\de_{\pm}(u) := \pm \left( {\rm e}^{\pm 2\pi \ii u} -1 \right) .
\end{equation}
By analogy with the case of flows,
we shall call \textit{Lindtsedt series} the series (\ref{eq:2.6}).

%%%%%%%%%%%%%%%%%%%%%%%%%%%%%%%%%%%%%%%%%%%%%%%%%
\begin{defi} \label{def:2.1}
We say that the formal series (\ref{eq:2.6}) is a \emph{formal solution}
to equation (\ref{eq:2.4}) if there exist well-defined coefficients $h^{(k)}_{\nu}$ 
that solve (\ref{eq:2.7}) to any order $k\in\NNN$.
\end{defi}
%%%%%%%%%%%%%%%%%%%%%%%%%%%%%%%%%%%%%%%%%%%%%%%%%

As (\ref{eq:2.7}) shows, to compute the coefficients $h^{(k)}_{\nu}$
one has to deal with infinitely many sums (over the Fourier labels).
Therefore (\ref{eq:2.7}) can be formally solved only if the sums can be performed
and the right hand side vanishes whenever $\de(\om\nu)=0$. In this section
we want to show that this occurs if we impose a mild Diophantine condition on $\om$.

We require for $\om$ to satisfy the Bryuno condition $\gotB(\om)<\io$, with $\gotB(\om)$
defined in (\ref{eq:1.4}). In particular this implies that
one has $\de(\om\nu)=0$ if and only if $\nu=0$. The following result holds.

%%%%%%%%%%%%%%%%%%%%%%%%%%%%%%%%%%%%%%%%%%%%%%%%%
\begin{prop} \label{prop:2.2}
Assume that $\om$ satisfies the Bryuno condition.
Then there exists a formal solution to (\ref{eq:2.4}) of the form (\ref{eq:2.6}).
The coefficients $h^{(k)}_{\nu}$ of the formal solution are uniquely defined
by requiring for the formal solution to have zero average, i.e. $h^{(k)}_{0}=0$ for all $k \ge 1$.
Moreover for any $\xi_{1},\xi_{2}\ge 0$, with $\xi_{1}+\xi_{2}<\xi$,
and all $k\in\NNN$ there exists a positive constant $C(k,\xi,\xi_{1},\xi_{2})$ such that
one has $|h^{(k)}_{\nu}| \le \xi_{2}^{-2k} C(k,\xi,\xi_1,\xi_2) \, {\rm e}^{-\xi_{1}|\nu|}$ 
for all $\nu\in\ZZZ$.
\end{prop}
%%%%%%%%%%%%%%%%%%%%%%%%%%%%%%%%%%%%%%%%%%%%%%%%%

We note since now that, in order to prove that a formal solution of
the form (\ref{eq:2.6}) exists, a condition even weaker than
the Bryuno condition  can be taken; see Remark \ref{rmk:2.7} below.
Of course Proposition \ref{prop:2.2} is not enough to prove the existence of the 
conjugation, because one still has to control that the constants $C(k,\xi,\xi_1,\xi_2)$
do not grow too fast in $k$. To check this a more detailed analysis is needed, 
as will be performed in Section \ref{sec:3}.

To prove Proposition \ref{prop:2.2} we shall use a tree expansion for the conjugation;
see also \cite{Ge1,GBG,Ge4} for an introduction to the subject.
A graph is a set of points and lines connecting them. A tree $\theta$ is
a graph with no cycle, such that all the lines are oriented toward a single point
(\emph{root}) which has only one incident line $\ell_{0}$ (\emph{root line}). 
All the points in a tree except the root are called \emph{nodes} or \emph{vertices}. 
The orientation of the lines in a tree induces a partial ordering relation ($\preceq$)
between the nodes and the lines: we can imagine that each line carries an arrow pointing
toward the root. Given two nodes $v$ and $w$, we shall write $w \prec v$ every time $v$
is along the oriented path of lines which connects $w$ to the root.
We denote by $V(\theta)$ and $L(\theta)$ the sets of nodes and lines in $\theta$,
respectively. Since a line $\ell\in L(\theta)$ is uniquely identified by the node $v$
which it leaves, we may write $\ell=\ell_{v}$ and say that $\ell$ exits $v$. We write
$\ell_{w} \prec \ell_{v}$ if $w\prec v$, and $w\prec \ell$ if $\ell=\ell_{v}$, with $w\preceq v$.
For any line $\ell\in L(\theta)$ the set of nodes $v' \prec \ell$ and the set of lines
$\ell'\preceq \ell$ form a tree $\theta'$: we shall say that $\theta'$ is a subtree
of $\theta$ with root line $\ell$; if $\ell$ enters the node $v$ we say that $\theta'$ enters $v$.

If $\ell$ and $\ell'$ are two comparable lines with $\ell'\prec \ell$, we
denote by $P(\ell,\ell')$ the oriented path of lines and nodes connecting $\ell'$ to $\ell$,
with $\ell$ and $\ell'$ not included (in particular $P(\ell,\ell')=\emptyset$
if $\ell'$ enters the node $\ell$ exits). If $w \prec v$ we denote by $\calP(v,w)$ 
the oriented path connecting $w$ to $v$, with $w$ and $v$ being included.

With each node $v\in V(\theta)$ we associate a \emph{mode} label 
$\nu_{v}\in\ZZZ$ and with each line $\ell\in L(\theta)$ 
we associate a \emph{momentum} label $\nu_{\ell}\in\ZZZ$, with the constraint
\begin{equation} \label{eq:2.9}
\nu_{\ell_{v}}= \sum_{w \preceq v} \nu_{v}  .
\end{equation}
With each line $\ell \prec \ell_{0}$ we associate also
a further label $\beta_{\ell}=\pm$, and we call $L_{v}^{\beta}(\theta)$ the set of lines $\ell$
entering the node $v\in V(\theta)$ such that $\beta_{\ell}=\beta$;
set $p_{v}=|L_{v}^{+}(\theta)|$ and $q_{v}=|L_{v}^{-}(\theta)|$.

We say that two trees are equivalent if they can be transformed
into each other by continuously deforming the lines in such a way that
these do not cross each other and all the labels match.
This provides an equivalence relation on the set of the trees.
From now on we shall call trees tout court such equivalence classes. 

We define the \emph{propagator} associated with the line $\ell$ as
\begin{equation} \label{eq:2.10}
G_{\ell} := \begin{cases}
\displaystyle{ \frac{1}{\de(\om\nu_{\ell})} } ,
& \nu_{\ell} \neq 0 , \\
1 & \nu_{\ell} = 0 ,
\end{cases}
\end{equation}
with $\de(u)$ defined in (\ref{eq:2.0}), and the \emph{node factor} associated with the node $v$ as
\begin{equation} \label{eq:2.11}
A_{v} :=  \frac{1}{p_v!q_v!} \Big( (\ii\nu_v)^{p_v+1} \partial_{z}^{q_v} 
+ \left( {\rm e}^{-2 \pi \ii \om\nu_{\ell_{v}}} - 1 \right) (\ii\nu_v)^{p_v} \partial_{z}^{q_v+1} \Big)
\s_{\n_{v}}(2\pi\om) \prod_{\ell\in L_{v}^{-}(\theta)} \de_{+}(\om\nu_{\ell}) ,
\end{equation}
where we have used that
\begin{equation} \label{eq:2.12}
\nu_{\ell_{v}} = \nu_{v} + \sum_{\ell\in L_{v}^{+}(\theta)} \nu_{\ell} +
\sum_{\ell\in L_{v}^{-}(\theta)} \nu_{\ell} . 
\end{equation}

Then we define the value of the tree $\theta$ as
\begin{equation} \label{eq:2.13}
\Val(\theta) = \Big( \prod_{v\in V(\theta)} A_{v} \Big)
\Big( \prod_{\ell\in L(\theta)} G_{\ell} \Big) 
\end{equation}
and call $\TT_{k,\nu}$ the set of all trees of order $k$ (that is with $k$ nodes)
and with momentum $\nu$ associated with the root line. Note that (\ref{eq:2.13})
is well-defined for $\om$ satisfying $\gotB(\om)<\io$.

%%%%%%%%%%%%%%%%%%%%%%%%%%%%%%%%%%%%%%%%%%%%%%%%%
\begin{lemma} \label{lem:2.3}
Set $h^{(k)}_{0}=0$ for all $k \ge 1$.
Assume that $\displaystyle{\sum_{\theta\in\TT_{k,0}}\Val(\theta)=0}$ $\forall k \le k_{0}$
for some $k_{0}\in\NNN$.

\vskip-.4truecm
\noindent Then, by setting
\begin{equation} \nonumber
h^{(k)}_{\nu} = \sum_{\theta\in \TT_{k,\nu}} \Val (\theta) , \quad \nu\neq 0 ,
\end{equation}
for all $k\le k_{0}$ and all $\nu\in\ZZZ$,
one obtains a solution up to order $k_{0}$ to (\ref{eq:2.7}).
\end{lemma}
%%%%%%%%%%%%%%%%%%%%%%%%%%%%%%%%%%%%%%%%%%%%%%%%%

%%%%%%%%%%%%%%%%%%%%%%%%%%%%%%%%%%%%%%%%%%%%%%%%%
\prova Equation (\ref{eq:2.7}) allows us to express the coefficients $h^{(k)}_{\nu}$, 
for $\nu\neq0$, in terms of coefficients $h^{(k')}_{\nu'}$, with $k'<k$ and $\nu'\neq0$
(since we are assuming $h^{(k)}_0=0$ for all $k\ge 1$).
By iterating, all the coefficients are eventually given by expressions
involving only coefficients $h^{(1)}_{\nu'}$, with $\nu'\neq 0$.
It is straightforward to see that such coefficients can be expressed as
\begin{equation} \label{eq:2.13b}
h^{(k)}_{\nu'} = \sum_{\theta\in \TT_{k,\nu'}} \Val (\theta) , \quad \nu'\neq 0 ,
\end{equation}
with $k=1$, provided $\Val(\theta)$ is defined as in (\ref{eq:2.13}). Then 
one easily checks by induction that $h^{(k)}_{\nu'}$ is given by (\ref{eq:2.13b})
for all $k\ge 1$ and all $\nu'\neq0$. To see that the definition of the coefficients
through (\ref{eq:2.13b}) makes sense, one uses that sum over the mode labels
can be performed thanks to (\ref{eq:2.2}) and the Bryuno condition $\gotB(\om)<\io$.

As a consequence, the equations (\ref{eq:2.7})
are satisfied for all $k\in\NNN$ and all $\nu\neq 0$.
The assumption ensures that, for all $k\le k_0$, the right hand
side of (\ref{eq:2.7}) vanishes for $\nu=0$.
Therefore the equations (\ref{eq:2.7}) are satisfied for all $k\le k_0$ and all $\nu\in\ZZZ$.
\qed

%%%%%%%%%%%%%%%%%%%%%%%%%%%%%%%%%%%%%%%%%%%%%%%%%
\begin{rmk} \label{rmk:2.4bis}
\emph{
The proof of Lemma \ref{lem:2.3} shows that, since any line $\ell$ of any tree $\theta$
can be considered as the root line of a subtree $\theta_{\ell}$, then, by construction,
$\Val(\theta_{\ell})$ is a contribution to $h^{(k_{\ell})}_{\nu_{\ell}}$, where $k_{\ell}$
is the order of $\theta_{\ell}$. In particular, the assumption $h^{(k)}_{0}=0$
implies that, for any tree $\theta$ and any line $\ell \in L(\theta)$, one has $\nu_{\ell}\neq 0$.
}
\end{rmk}
%%%%%%%%%%%%%%%%%%%%%%%%%%%%%%%%%%%%%%%%%%%%%%%%%

%%%%%%%%%%%%%%%%%%%%%%%%%%%%%%%%%%%%%%%%%%%%%%%%%
\begin{lemma} \label{lem:2.4}
One has $\displaystyle{\sum_{\theta\in\TT_{k,0}}\Val(\theta)=0}$ for all $k\ge 1$.
\end{lemma}
%%%%%%%%%%%%%%%%%%%%%%%%%%%%%%%%%%%%%%%%%%%%%%%%%

The proof of Lemma \ref{lem:2.4} is given in Appendix \ref{app:a}.

%%%%%%%%%%%%%%%%%%%%%%%%%%%%%%%%%%%%%%%%%%%%%%%%%
\begin{rmk} \label{rmk:2.5}
\emph{
Lemma \ref{lem:2.4}, together with Lemma \ref{lem:2.3},
implies that one can define the coefficients
$h^{(k)}_{\nu}$ as in Lemma \ref{lem:2.3} for all $k\in\NNN$.
Moreover, by Remark \ref{rmk:2.4bis}, if $\theta\in\TT_{k,\nu}$
and $\ell_{0}$ is the root line of $\theta$,
then $\nu_{\ell} \neq 0$ for all $\ell \in L(\theta) \setminus\{\ell_{0}\}$. 
}
\end{rmk}
%%%%%%%%%%%%%%%%%%%%%%%%%%%%%%%%%%%%%%%%%%%%%%%%%

The following result yields immediately Proposition \ref{prop:2.2}.

%%%%%%%%%%%%%%%%%%%%%%%%%%%%%%%%%%%%%%%%%%%%%%%%%
\begin{lemma} \label{lem:2.6}

For any $\xi_{1},\xi_{2}\ge 0$, with $\xi_{1}+\xi_{2}<\xi$, and any $k\ge 1$
there exists a constant $C_{0}(k,\xi,\xi_{1},\xi_{2})$ such that one has
\begin{equation} \nonumber
\sum_{\theta\in\TT_{k,\nu}}
\left| \Val(\theta) \right| \le \xi_{2}^{-2k} 
C_{0}(k,\xi,\xi_{1},\xi_{2}) \, {\rm e}^{-\xi_{1}|\nu|}
\end{equation}
for any $\nu\in\ZZZ\setminus\{0\}$ .
\end{lemma}
%%%%%%%%%%%%%%%%%%%%%%%%%%%%%%%%%%%%%%%%%%%%%%%%%

%%%%%%%%%%%%%%%%%%%%%%%%%%%%%%%%%%%%%%%%%%%%%%%%%
\prova 
With the notations in (\ref{eq:2.2}), for any $\xi_{2}\in(0,\xi)$
each node factor $A_{v}$ can be bounded by $3\,\Xi\,(2\,\xi_{2}^{-1})^{(p_{v}+1)}
2^{q_{v}}\rho^{-(q_{v}+1)}{\rm e}^{-(\xi-\xi_{2})|\nu_{v}|}$. By writing $\xi-\xi_{2}=
\xi_{1}+(\xi-\xi_{1}-\xi_{2})$ in ${\rm e}^{-(\xi-\xi_{2})|\nu_{v}|}$,
one can extract a factor ${\rm e}^{-\xi_{1}|\nu_{v}|}$ per each node.
Moreover one has
\begin{equation} \nonumber
\sum_{v\in V(\theta)} \left( p_{v}+q_{v} \right) = k-1 .
\end{equation}
Then one uses the Bryuno condition in order to bound the sum over
all the mode labels of the product of factors ${\rm e}^{-(\xi-\xi_{1}-\xi_{2})
|\nu_{v}|}$ times the product of propagators by a
constant $C_{1}(k,\xi,\xi_{1},\xi_{2})$; one reasons as in Appendix H in \cite{CG},
which we refer to for more details.
Finally one has to sum over all the other labels and this produces a further
factor $C_{2}^{k}$ for a suitable constant $C_{2}$. Then the assertion follows with
$C_{0}(k,\xi,\xi_{1},\xi_{2})=(24\,\Xi\,C_{2}\rho^{-2})^{k}C_{1}(k,\xi,\xi_{1},\xi_{2})$.
\qed
%%%%%%%%%%%%%%%%%%%%%%%%%%%%%%%%%%%%%%%%%%%%%%%%%

%%%%%%%%%%%%%%%%%%%%%%%%%%%%%%%%%%%%%%%%%%%%%%%%%
\begin{rmk} \label{rmk:2.7}
\emph{
As a matter of fact, the existence of the formal solution (\ref{eq:2.6})
follows by only assuming that $q_{n}^{-1}\log q_{n+1}\to 0$ as $n\to\io$;
see also \cite{CG} for a similar analysis.
}
\end{rmk}
%%%%%%%%%%%%%%%%%%%%%%%%%%%%%%%%%%%%%%%%%%%%%%%%%

\vspace{.3truecm}

%%%%%%%%%%%%%%%%%%%%%%%%%%%%%%%%%%%%%%%%%%%%%%%%%
\begin{rmk} \label{rmk:2.8}
\emph{
The formal solution to (\ref{eq:2.4}) found through the construction above has
been obtained by arbitrarily imposing that $h^{(k)}_{0}=0$ for all $k\in\NNN$.
Therefore we have no uniqueness of the solution: the solution is unique
if we require for its average to vanish. In principle one could fix arbitrarily the constants
$h^{(k)}_0$ and the existence of the solution could still be proved by reasoning as above: 
this would simply shift the origin of the parametrisation of the invariant curve. 
However the choice we made is a natural one and simplifies the analysis.
}
\end{rmk}
%%%%%%%%%%%%%%%%%%%%%%%%%%%%%%%%%%%%%%%%%%%%%%%%%

%%%%%%%%%%%%%%%%%%%%%%%%%%%%%%%%%%%%%%%%%%%%%%%%%
%%%%%%%%%%%%%%%%%%%%%%%%%%%%%%%%%%%%%%%%%%%%%%%%%
\zerarcounters 
\section{Analyticity of the conjugation: multiscale analysis}
\label{sec:3} 
%%%%%%%%%%%%%%%%%%%%%%%%%%%%%%%%%%%%%%%%%%%%%%%%%
%%%%%%%%%%%%%%%%%%%%%%%%%%%%%%%%%%%%%%%%%%%%%%%%%

To prove convergence of the power series in $\e$ for the conjugation
we need a more careful analysis, which requires multiscale techniques.

We introduce a $C^{\io}$ partition of unity as follows. 
For $b>a>0$ let $\chi_{a,b}$ be a $C^{\io}$ non-increasing function
defined on $(0,+\io)$ such that
\begin{equation} \label{eq:3.1}
\chi_{a,b}(u) = \begin{cases}
1 , & u \le a , \\
0 , & u \ge b .
\end{cases}
\end{equation}
For $n\ge 1$ define
\begin{equation} \label{eq:3.2}
\begin{split}
& \delta_{n}^{-} := \frac{1}{4} \min\left\{ 
\frac{q_{n+1}-q_{n}}{32\,q_{n}q_{n+1}} , \frac{1}{32q_{n}} \right\} ,
\quad \delta_{n}^{+} := \frac{1}{4} \min\left\{ 
\frac{q_{n}-q_{n-1}}{32\,q_{n-1}q_{n}} , \frac{1}{32q_{n}} \right\} , \\
& a_{n} := \frac{1}{32q_{n}}-\de_{n}^{-} , \qquad
\quad b_{n} := \frac{1}{32q_{n}}+\de_{n}^{+} .
\end{split}
\end{equation}
and set $\chi_{n}=\chi_{a_{n},b_{n}}$. Then we define
\begin{equation} \label{eq:3.3}
\Psi_{0}(u)=1-\chi_{1}(u) ,  \qquad
\Psi_{n}(u) = \chi_{n}(u)-\chi_{n+1}(u) , \quad n \ge 1 .
\end{equation}
%

%%%%%%%%%%%%%%%%%%%%%%%%%%%%%%%%%%%%%%%%%%%%%%%%%
\begin{rmk} \label{rmk:3.1}
\emph{
Consider the partition of unity of $(0,+\io)$ consisting of the characteristic functions of
the intervals $(1/32q_{n+1},1/32q_{n}]$, $n\ge 0$, and of the half-line $(1/32q_{1},+\io)$.
Then $\{\Psi_{n}\}_{n\ge 0}$ is a smoothed version of such a partition.
}
\end{rmk}
%%%%%%%%%%%%%%%%%%%%%%%%%%%%%%%%%%%%%%%%%%%%%%%%%

We associate with each line $\ell\in L(\theta)$ a further label $n_{\ell}\in\ZZZ_{+}$,
called the \emph{scale label},
and redefine the propagator of the line $\ell$ as
\begin{equation} \label{eq:3.4}
\calG_{\ell} = \Psi_{n_{\ell}}(\|\om\nu_{\ell}\|) \, G_{\ell} =
\frac{\Psi_{n_{\ell}}(\|\om\nu_{\ell}\|)}{\de(\om\nu_{\ell})} ,
\end{equation}
where
\begin{equation} \label{eq:3.5}
\|u\| :=\inf_{\nu\in\ZZZ} \left| u - \nu \right| .
\end{equation}
%

%%%%%%%%%%%%%%%%%%%%%%%%%%%%%%%%%%%%%%%%%%%%%%%%%
\begin{lemma} \label{lem:3.2}
If $\calG_{\ell} \neq 0$, then $1/64q_{n_{\ell}+1} \!\le\! \|\om\nu_{\ell}\|
\!\le\! 1/6q_{n_{\ell}}$ if $n_{\ell} \!\ge\! 1$ and $1/64q_{1} \!\le\! \|\om\nu_{\ell}\|$
if $n_{\ell} \!=\! 0$.
\end{lemma}
%%%%%%%%%%%%%%%%%%%%%%%%%%%%%%%%%%%%%%%%%%%%%%%%%

%%%%%%%%%%%%%%%%%%%%%%%%%%%%%%%%%%%%%%%%%%%%%%%%%
\prova By definition of the functions $\Psi_{n}$.
\qed
%%%%%%%%%%%%%%%%%%%%%%%%%%%%%%%%%%%%%%%%%%%%%%%%%

%%%%%%%%%%%%%%%%%%%%%%%%%%%%%%%%%%%%%%%%%%%%%%%%%
\begin{defi} \label{def:3.3}
For any $\ell\in L(\theta)$ we call
\begin{equation} \nonumber
\zeta_{\ell} = \min\{ n \in \ZZZ_{+} : \Psi_{n}(\|\om\nu_{\ell}\|) \neq 0 \}
\end{equation}
the \emph{minimum scale} of $\ell$.
\end{defi}
%%%%%%%%%%%%%%%%%%%%%%%%%%%%%%%%%%%%%%%%%%%%%%%%%

%%%%%%%%%%%%%%%%%%%%%%%%%%%%%%%%%%%%%%%%%%%%%%%%%
\begin{rmk} \label{rmk:3.4}
\emph{
For any line $\ell \in L(\theta)$ such that $\calG_{\ell}\neq0$ one has
either $n_{\ell}=\zeta_{\ell}$ or $n_{\ell}=\zeta_{\ell}+1$.
}
\end{rmk}
%%%%%%%%%%%%%%%%%%%%%%%%%%%%%%%%%%%%%%%%%%%%%%%%%

We denote by $\Theta_{k,\nu}$ the set of all labelled trees with such extra labels
$n_{\ell},\zeta_{\ell}$ for all $\ell\in L(\theta)$,
with $k$ nodes and momentum $\nu$ associated with the root line.
Given a tree $\theta\in\Theta_{k,\nu}$ redefine the value of $\theta$ as
\begin{equation} \label{eq:3.6}
\Val(\theta) = \Big( \prod_{v\in V(\theta)} A_{v} \Big)
\Big( \prod_{\ell\in L(\theta)} \calG_{\ell} \Big) .
\end{equation}
%

%%%%%%%%%%%%%%%%%%%%%%%%%%%%%%%%%%%%%%%%%%%%%%%%%
\begin{defi} \label{def:3.5}
We say that a line $\ell$ with scale $n_{\ell}=n$ satisfies the \emph{support property} if
$1/128q_{n+1} \le \|\om\nu_{\ell}\|\le 1/8q_{n}$.
\end{defi}
%%%%%%%%%%%%%%%%%%%%%%%%%%%%%%%%%%%%%%%%%%%%%%%%%

%%%%%%%%%%%%%%%%%%%%%%%%%%%%%%%%%%%%%%%%%%%%%%%%%
\begin{rmk} \label{rmk:3.7}
\emph{
Lemma \ref{lem:3.2} implies that if $\Val(\theta) \neq 0$ then all lines
$\ell\in L(\theta)$ satisfy the support property. However, as it will be shown
later, it is convenient to consider also trees with zero value,
but still satisfying the support property. This will be due to the fact that
in the renormalisation procedure described in Section \ref{sec:4} the momenta
$\nu_{\ell}$ will be replaced by new momenta, to be denoted $\nu_{\ell}(\bt)$,
so that the scales $n_{\ell}$ will be determined by these new momenta
through the conditions $\Psi_{n_{\ell}}(\|\om\nu_{\ell}(\bt)\|) \neq 0$.
Nevertheless the old momenta $\nu_{\ell}$ will be still found to
satisfy the support property in terms of the scales $n_{\ell}$.
}
\end{rmk}
%%%%%%%%%%%%%%%%%%%%%%%%%%%%%%%%%%%%%%%%%%%%%%%%%

The following result is proved as Lemma \ref{lem:2.3} and Lemma \ref{lem:2.4}.

%%%%%%%%%%%%%%%%%%%%%%%%%%%%%%%%%%%%%%%%%%%%%%%%%
\begin{lemma} \label{lem:3.6}
The coeffiecients
\begin{equation} \nonumber
h^{(k)}_{0}=0 \qquad \hbox{and} \qquad
h^{(k)}_{\nu} = \sum_{\theta\in \Theta_{k,\nu}} \Val (\theta) , \quad \nu \neq 0 ,
\end{equation}
for all $k\ge 1$ and all $\nu\in\ZZZ$, formally solve (\ref{eq:2.7}) to all orders $k\in\NNN$.
\end{lemma}
%%%%%%%%%%%%%%%%%%%%%%%%%%%%%%%%%%%%%%%%%%%%%%%%%

%%%%%%%%%%%%%%%%%%%%%%%%%%%%%%%%%%%%%%%%%%%%%%%%%
\begin{defi} \label{def:3.8}
We define a \emph{cluster} $T$ on scale $n$ as a maximal connected set
of nodes and lines connecting them such that all lines have scale $\le n$
and at least one of them has scale $n$. 
\end{defi}
%%%%%%%%%%%%%%%%%%%%%%%%%%%%%%%%%%%%%%%%%%%%%%%%%

Given any subgraph $T$ of a tree $\theta$, denote by $L(T)$ and $V(T)$
the set of lines and the set of nodes, respectively, in $T$.
If $T$ is a cluster, $L(T)$ will be called the set of \emph{internal lines} of $T$,
while the lines which connect one node $v\in V(T)$ to a node $w\notin V(T)$
will be called the \emph{external lines} of $T$: they will be said to enter $T$
if the line is oriented toward to the node belonging to $V(T)$ and to exit $T$ otherwise.
A cluster can have only either one or no exiting line.

If $T$ has one exiting line $\ell_{T}$ and only one entering line $\ell_{T}'$
call $\calP_{T}$ the oriented path of lines and nodes
connecting $\ell_{T}'$ to $\ell_{T}$ and define $n_{T}:=\min\{n_{\ell_{T}},n_{\ell_{T}'}\}$.
If $T$ is a subgraph of $\theta$, consisting of all nodes and lines
preceding the line $\ell_{1}$  but not the line $\ell_{2}$, then we shall still
denote by $\calP_{T}$ the oriented path connecting $\ell_{2}$ to $\ell_{1}$.
Finally, for $\ell\in\calP_{T}$, call $w(\ell)$ and $w'(\ell)$ the nodes along $\calP_{T}$
which $\ell$ enters and exits, respectively.

%%%%%%%%%%%%%%%%%%%%%%%%%%%%%%%%%%%%%%%%%%%%%%%%%
\begin{rmk} \label{rmk:3.9}
\emph{
If $T$ is a cluster on scale $n$ with one entering line and one exiting line,
one has $n_{T} \ge n+1$.
}
\end{rmk}
%%%%%%%%%%%%%%%%%%%%%%%%%%%%%%%%%%%%%%%%%%%%%%%%%

Set also
\begin{equation} \label{eq:3.7}
M(\theta) = \sum_{v\in V(\theta)} |\nu_{v}| , \qquad
M(T) = \sum_{v\in V(T)} |\nu_{v}| ,
\end{equation}
where $T$ can be any subgraph of $\theta$.

%%%%%%%%%%%%%%%%%%%%%%%%%%%%%%%%%%%%%%%%%%%%%%%%%
\begin{defi} \label{def:3.10}
A cluster $T$ on scale $n$ will be called a \emph{self-energy cluster} if\\
(1) it has one exiting line $\ell_{T}$ and only one entering line $\ell_{T}'$,\\
(2) $\nu_{\ell_{T}'}=\nu_{\ell_{T}}$,\\
(3) $\nu_{\ell} \neq \nu_{\ell_{T}}$ for all $\ell\in\calP_{T}$,\\
(4) $M(T) < q_{n_{T}}$.
\end{defi}
%%%%%%%%%%%%%%%%%%%%%%%%%%%%%%%%%%%%%%%%%%%%%%%%%

%%%%%%%%%%%%%%%%%%%%%%%%%%%%%%%%%%%%%%%%%%%%%%%%%
\begin{rmk} \label{rmk:3.11}
\emph{
The condition (2) in Definition \ref{def:3.10} implies that the sum of the mode labels
of the nodes of a self-energy cluster $T$ is zero, that is
\begin{equation} \nonumber
\sum_{v\in V(T)} \nu_{v} =0 .
\end{equation}
}
\end{rmk}
%%%%%%%%%%%%%%%%%%%%%%%%%%%%%%%%%%%%%%%%%%%%%%%%%

We say that $\ell\in L(\theta)$ is a \emph{resonant line} if $\ell$ exits a self-energy cluster.
Let us denote by $N_{n}^{\bullet}(\theta)$ the number of
non-resonant lines $\ell\in L(\theta)$ with $\zeta_{\ell} = n$
and by $P_{n}(\theta)$ the number of self-energy clusters on scale $n$ in $\theta$.

%%%%%%%%%%%%%%%%%%%%%%%%%%%%%%%%%%%%%%%%%%%%%%%%%
\begin{lemma} \label{lem:3.12}
Given $\nu\in\ZZZ$, if $\|\om\nu\|\le 1/4q_{n}$ then either $\nu=0$ or $|\nu|\ge q_{n}$.
\end{lemma}
%%%%%%%%%%%%%%%%%%%%%%%%%%%%%%%%%%%%%%%%%%%%%%%%%

The proof of Lemma \ref{lem:3.12} can be found in \cite{D}.

%%%%%%%%%%%%%%%%%%%%%%%%%%%%%%%%%%%%%%%%%%%%%%%%%
\begin{lemma} \label{lem:3.13}
Assume that all lines $\ell\in L(\theta)$ satisfy the support property. Then one has
\begin{equation} \nonumber
K_{n}(\theta):= N^{\bullet}_{n}(\theta) + P_{n}(\theta) \le \frac{2 M(\theta)}{q_{n}} .
\end{equation}
\end{lemma}
%%%%%%%%%%%%%%%%%%%%%%%%%%%%%%%%%%%%%%%%%%%%%%%%%

%%%%%%%%%%%%%%%%%%%%%%%%%%%%%%%%%%%%%%%%%%%%%%%%%
\proof We prove by induction that for any tree $\theta$ one has
\begin{subequations} \label{eq:3.8}
\begin{align}
& K_{n}(\theta) = 0 , \hskip2.65truecm  M(\theta) < q_{n} ,
\label{eq:3.8a} \\
& K_{n}(\theta) \le \frac{2M(\theta)}{q_{n}} -1 , \hskip1.truecm M(\theta) \ge q_{n} ,
\label{eq:3.8b}
\end{align}
\end{subequations}
First of all note that if $M(\theta)<q_{n}$ one has $|\nu_{\ell}|<q_{n}$ and hence,
by Lemma \ref{lem:3.12}, $n_{\ell}<n$ for all $\ell \in L(\theta)$: 
therefore $K_{n}(\theta) \ge 1$ requires $M(\theta)\ge q_{n}$.
In particular this proves (\ref{eq:3.8a}).
To prove (\ref{eq:3.8b}) we proceed by induction on the order $k$ of the tree.
For $k=1$ the bound holds by the previous argument. Let assume that (\ref{eq:3.8b})
holds for any tree of order $k<k_{0}$ and consider a tree of order $k_{0}$.
Let us denote by $\ell_{0}$ the root line of $\theta$ and call
$\ell_{1},\ldots,\ell_{p}$ the lines with minimum scale $\ge n$ (if any)
which are closest to $\ell_{0}$ and $\theta_{1},\ldots,\theta_{p}$
the subtrees with root lines $\ell_{1},\ldots,\ell_{p}$, respectively.
Note that $n_{\ell_{0}}\ge \zeta_{\ell_{0}}$, so $\ell_{0}$
can exit a self-energy cluster on scale $n$ only if $\zeta_{\ell_{0}}\ge n$
and $n_{\ell_{0}}\ge n+1$.

\begin{enumerate}

\item If $\zeta_{\ell_{0}} \neq n$ and $\ell_{0}$ does not exit a self-energy cluster
on scale $n$, then one has $K_{n}(\theta)=0$ if $p=0$ and
$K_{n}(\theta)=K_{n}(\theta_{1})+\ldots+K_{n}(\theta_{p})$ if $p\ge 1$.
In the first case the bound (\ref{eq:3.8b}) trivially holds, while in the
second case it follows from the inductive hypothesis.

\item If $\zeta_{\ell_{0}} = n$ and $\ell_{0}$ is non-resonant, then one has
$K_{n}(\theta)=1+K_{n}(\theta_{1})+\ldots+K_{n}(\theta_{p})$. If $p=0$
the bound is trivially satisfied and if $p\ge 2$ it follows immediately
from the inductive hypothesis. If $p=1$ let $T$ be the subgraph formed
by all lines and nodes of $\theta$ which precede $\ell_{0}$ but not $\ell_{1}$.
If $T$ is not a cluster then it must contain at least one line $\ell$
with $n_{\ell}=n$ and $\zeta_{\ell}=n-1$. Then $\|\om\nu_{\ell}\|\le 1/8q_{n}$.
If $\ell$ is not along the path $\calP_{T}$
connecting $\ell_{1}$ to $\ell_{0}$ then $M(T)\ge |\nu_{\ell}|\ge q_{n}$;
if $\ell$ is along the path $\calP_{T}$ then $\nu_{\ell}\neq\nu_{\ell_{1}}$
(because $\zeta_{\ell}=n-1$, while $\zeta_{\ell_{1}}\ge n$),
so that  $\|\om(\nu_{\ell}-\nu_{\ell_{1}})\| \le
\| \om \nu_{\ell} \| + \| \om \nu_{\ell_{1}} \| \le 1/4q_{n}$ implies
$|\nu_{\ell}-\nu_{\ell_{1}}|\ge q_{n}$ and hence $M(T)\ge q_{n}$.
If $T$ is a cluster then either $M(T)\ge q_{n}$ or
$\nu_{\ell_{0}}\neq\nu_{\ell_{1}}$ (as $\ell_{0}$ is non-resonant,
$T$ cannot be a self-energy cluster and hence $M(T)<q_{n}$ would require
$\nu_{\ell_{0}}\neq \nu_{\ell_{1}}$,
because otherwise there should be at least one line $\ell\in\calP_{T}$
such that $\nu_{\ell}=\nu_{\ell_{1}}$, which is not possible
since $\zeta_{\ell}\le n-1$ and $\zeta_{\ell_{1}} =\zeta_{\ell_{0}}=n$).
On the other hand $\nu_{\ell_{0}}\neq \nu_{\ell_{1}}$ would imply
$\| \om \left( \nu_{\ell_{0}} - \nu_{\ell_{1}} \right) \| \le
\| \om \nu_{\ell_{0}} \| + \| \om \nu_{\ell_{1}} \|  \le 1/4q_{n}$ and hence
$M(T) \ge q_{n}$ once more.
Therefore, in all cases one has $M(\theta)-M(\theta_{1}) =
M(T) \ge q_{n}$, so that $K_{n}(\theta)=1+K(\theta_{1}) \le
2M(\theta_{1})/q_{n} \le 2M(\theta)/q_{n}-1$
and the bound (\ref{eq:3.8b}) follows.

\item If $\zeta_{\ell_{0}}=n$ and $\ell_{0}$ exits a self-energy cluster $T$
on scale $\neq n$, then one has
$K_{n}(\theta)=K_{n}(\theta_{1})+\ldots+K_{n}(\theta_{p})$ and once more
(\ref{eq:3.8b}) is trivial if $p=0$ and
follows from the inductive hypothesis if $p\ge 1$.

\item If $\ell_{0}$ exits a self-energy cluster $T$ on scale $n$,
then $p\ge 1$ and $K_{n}(\theta)=1+K(\theta_{1})+\ldots+K(\theta_{p})$.
If $p\ge2$ the bound follows by the inductive hypothesis.
If $p=1$ then either $\ell_{1}=\ell_{T}'$ or $\ell_{1} \in \calP_{T}$
(because $\ell_{T}'$ has $n_{\ell_{T}'}\ge n+1$ and hence $\zeta_{\ell_{T}'}\ge n$).
If $\ell_{1}=\ell_{T}'$ then $T$ must contain at least one line $\ell$ such that
$n_{\ell}=n$ and $\zeta_{\ell}=n-1$: if $\ell\not\in\calP_{T}$
then $M(T) \ge |\nu_{\ell}| \ge q_{n}$, while if $\ell\in\calP_{T}$ then
either $M(T) \ge q_{n}$ or $\nu_{\ell}\neq\nu_{\ell_{1}}$ and hence
again $M(T)\ge |\nu_{\ell_{1}}-\nu_{\ell}| \ge q_{n}$.
If $\ell_{1} \in \calP_{T}$ then $n_{\ell_{1}}=n$ and either
$M(T) \ge q_{n}$ or $\nu_{\ell_{1}} \neq \nu_{\ell_{T}'}$: the latter
case yields again  $M(T) \ge |\nu_{\ell_{1}}-\nu_{\ell_{T}'}| \ge q_{n}$.
In all cases one finds $K_{n}(\theta) \le
2M_{n}(\theta_{1})/q_{n} \le 2M_{n}(\theta)/q_{n}-1$, so that (\ref{eq:3.8b}) follows.
\qed
\end{enumerate}
%%%%%%%%%%%%%%%%%%%%%%%%%%%%%%%%%%%%%%%%%%%%%%%%%

As we shall see in Section \ref{sec:4}, when bounding a tree value (\ref{eq:3.6}),
the product of node factors is easily controlled, while Lemma \ref{lem:3.13} guarantees
that the product of propagators would also be bounded proportionally to a constant to the
power $k$, if only the resonant lines could be neglected. In other words,
the small divisors can accumulate  only in the presence of self-energy clusters.
Then, we have to make sure that when this happens one has compensations
between the tree values. We end this section by providing the basic algebraic
cancellations which underlie such compensations; the proof is an extension
of the very argument used in the proof of Lemma \ref{lem:2.4} and is
easily described in terms of operations on trees.
Finally, in Section \ref{sec:4} we will show show how to implement iteratively the
cancellations in order to control the product of all propagators, including those of the
resonant lines, and hence complete the proof of convergence of the series (\ref{eq:2.6}).

We denote by $\gotS_{k,n}$ the set of all self-energy clusters $T$ on scale $n$ of order
$k$, that is with $k$ nodes. Given a self-energy cluster $T$ define the value of $T$ as
\begin{equation} \label{eq:3.9}
\Val_{T}(\om\nu_{\ell_{T}'}) =
\Big( \prod_{v\in V(T)} A_{v} \Big)
\Big( \prod_{\ell\in L(T)} \calG_{\ell} \Big) .
\end{equation}
%

%%%%%%%%%%%%%%%%%%%%%%%%%%%%%%%%%%%%%%%%%%%%%%%%%
\begin{rmk} \label{rmk:3.14}
\emph{
Given a self-energy cluster $T$, if we set $\nu=\nu_{\ell_{T}'}$ and define
\begin{equation} \label{eq:3.10}
\nu_{\ell}^{0} = \sum_{\substack{v \in V(T) \\ v \prec \ell}} \nu_{v} ,
\end{equation}
we have $\nu_{\ell}=\nu_{\ell}^{0}+\nu$ if $\ell\in\calP_{T}$
and $\nu_{\ell}=\nu_{\ell}^{0}$ if $\ell\notin\calP_{T}$.
Therefore $\Val_{T}(\om\nu)$ depends on $\om\nu$ only through the propagators and
the node factors of the lines and of the nodes, respectively, along the path $\calP_{T}$.
}
\end{rmk}
%%%%%%%%%%%%%%%%%%%%%%%%%%%%%%%%%%%%%%%%%%%%%%%%%

A self-energy cluster $T$ can contain other self-energy clusters:
call $\mathring{T}$ the set of lines and nodes in $T$ which are
outside any self-energy clusters contained inside $T$, and
denote by $L(\mathring{T})$ and $V(\mathring{T})$ the set of
lines and of nodes, respectively, in $\mathring{T}$. 

%%%%%%%%%%%%%%%%%%%%%%%%%%%%%%%%%%%%%%%%%%%%%%%%%
\begin{rmk} \label{rmk:3.15}
\emph{
By definition
of self-energy cluster one has
\begin{equation} \nonumber
\sum_{v\in V(\mathring{T})} \nu_{v} = 0
\end{equation}
for any self-energy cluster $T$.
}
\end{rmk}
%%%%%%%%%%%%%%%%%%%%%%%%%%%%%%%%%%%%%%%%%%%%%%%%%

Given a self-energy cluster $T$, call $\gotF(T)$ the set of all self-energy clusters
$T'\in\gotS_{k,n}$, where $k$ and $n$ are the order and scale, respectively, of $T$,
with $\nu_{\ell_{T'}'}=\nu_{\ell_{T}'}$.

%%%%%%%%%%%%%%%%%%%%%%%%%%%%%%%%%%%%%%%%%%%%%%%%%
\begin{lemma} \label{lem:3.16}
For any self-energy cluster $T$ one has
\begin{equation} \nonumber
\sum_{T'\in \gotF(T)} \Val_{T'}(0) = 0 , \qquad
\sum_{T'\in \gotF(T)} \partial_{u} \Val_{T'}(0) = 0 .
\end{equation}
\end{lemma}
%%%%%%%%%%%%%%%%%%%%%%%%%%%%%%%%%%%%%%%%%%%%%%%%%

The proof is in Appendix \ref{app:a}.

%%%%%%%%%%%%%%%%%%%%%%%%%%%%%%%%%%%%%%%%%%%%%%%%%
%%%%%%%%%%%%%%%%%%%%%%%%%%%%%%%%%%%%%%%%%%%%%%%%%
\zerarcounters 
\section{Analyticity of the conjugation: renormalisation procedure}
\label{sec:4} 
%%%%%%%%%%%%%%%%%%%%%%%%%%%%%%%%%%%%%%%%%%%%%%%%%
%%%%%%%%%%%%%%%%%%%%%%%%%%%%%%%%%%%%%%%%%%%%%%%%%

The modified tree expansion envisaged in Section \ref{sec:3}
will allow us to refine the bounds on the coefficients of Proposition \ref{prop:2.2} 
into the following result.

%%%%%%%%%%%%%%%%%%%%%%%%%%%%%%%%%%%%%%%%%%%%%%%%%
\begin{prop} \label{prop:4.1}
Assume that $\om$ satisfies the Bryuno condition.
Then there exists a solution $h(\psi)$ to (\ref{eq:2.4})
of the form (\ref{eq:2.6}). For any $\xi_{1},\xi_{2},\xi_{3}\ge 0$, with $\xi_{1}+\xi_{2}+\xi_{3}<\xi$,
there exists a positive constant $C_{0}(\xi_{3})$ such that for all $k\ge 1$
and $\nu\in\ZZZ$ the coefficients $h^{(k)}_{\nu}$ of the solution $h(\psi)$
satisfy the bounds $|h^{(k)}_{\nu}| \le (\xi-\xi_{1}-\xi_{2}-\xi_{3})^{-k}\xi_{2}^{-2k}
C_{0}^{k}(\xi_{3}) {\rm e}^{-\xi_{1}|\nu|}$ and hence the series (\ref{eq:2.6})
converges absolutely for $\e$ small enough.
\end{prop}
%%%%%%%%%%%%%%%%%%%%%%%%%%%%%%%%%%%%%%%%%%%%%%%%%

The rest of this section is devoted to proving Proposition \ref{prop:4.1},
which in turn yields immediately Theorem \ref{thm:2.0}.

We define $\gotT(\theta)$ as the set of self-energy clusters in $\theta$ and
$\gotT_{1}(\theta)$ as the set of maximal self-energy clusters $T$
contained in $\theta$, i.e. such that there is no self-energy cluster in $\theta$ containing $T$. 
Generally, for any self-energy cluster $T'$ in $\theta$, we denote by $\gotT_{1}(T')$
the set of maximal self-energy clusters (strictly) contained in $T'$.
Given a line $\ell\in L(\theta)$, either there is no self-energy cluster containing $\ell$
or there exist $p=p(\ell)\ge1$ self-energy clusters $T_{1},\ldots,T_{p}$ such that 
$T_{p}$ is the minimal self-energy cluster containing $\ell$ and
$T_{1} \supset T_{2} \supset \ldots \supset T_{p}$, with
$T_{j}\in\gotT_{1}(T_{j-1})$ for $j=2,\ldots,p$ and $T_{1} \in \gotT_{1}(\theta)$.
We call $\gotC_{\ell}(\theta):=\{T_{1},\ldots,T_{p}\}$
the \emph{cloud} of $\ell$ in $\theta$.

Define $\mathring{\theta}$ as the set of nodes and lines in $\theta$ which are outside
any self-energy cluster $T\in\gotT_{1}(\theta)$ and set
\begin{equation} \label{eq:4.1}
\Val(\mathring{\theta}) = 
\Big( \prod_{v\in V(\mathring{\theta})} A_{v} \Big)
\Big( \prod_{\ell\in L(\mathring{\theta})} \calG_{\ell} \Big) .
\end{equation}
We can write $\Val(\theta)$ in (\ref{eq:3.6}) as
\begin{equation} \label{eq:4.2}
\Val(\theta) = \Val(\mathring{\theta}) \prod_{T\in\gotT_{1}(\theta)} \Val_{T}(\om\nu_{\ell_{T}'}) .
\end{equation}
Define the \emph{localised value} of a self-energy cluster $T\in\gotT_{1}(\theta)$ as
\begin{equation} \label{eq:4.3}
\LL \Val_{T}(u) = 
\Val_{T}(0) + u \, \partial_{u} \Val_{T}(0)
\end{equation}
and the \emph{regularised value} of $T$ as $\RR\Val_{T}(u)=\Val_{T}(u)-\LL\Val_{T}(u)$.
One can write
\begin{equation} \label{eq:4.4}
\RR\Val_{T}(u) = u^{2} \int_{0}^{1} {\rm d}t_{T} \left( 1 - t_{T} \right)
\partial_{u}^{2} \Val_{T}(t_{T}u) ,
\end{equation}
where $t_{T}\in[0,1]$ will be called the \emph{interpolation parameter} associated with
the self-energy cluster $T$.
Then we associate with any $T\in\gotT_{1}(\theta)$ a label $\de_{T}\in\{\LL,\RR\}$ and
rewrite (\ref{eq:4.2}) as a sum of contributions
\begin{equation} \label{eq:4.5}
\Val(\mathring{\theta})
\Big( \prod_{\substack{T\in\gotT_{1}(\theta) \\ \de_{T}=\RR}}
\RR \Val_{T}(\om\nu_{\ell_{T}'}) \Big)
\Big( \prod_{\substack{T\in\gotT_{1}(\theta) \\ \de_{T}=\LL}}
\LL \Val_{T}(\om\nu_{\ell_{T}'}) \Big) .
\end{equation}

Recall the definition of the nodes $w(\ell)$ and $w'(\ell)$ before
Remark \ref{rmk:3.9}. Call $A_{w(\ell)}^{1}$ the node factor obtained
from $A_{w(\ell)}$ by replacing the factor $\de_{+}(\om\nu_{\ell})$
with $2\pi\ii \, {\rm e}^{2\pi\ii\om\nu_{\ell}}$ if $\be_{\ell}=-$ and
with $0$ if $\be_{\ell}=+$, and call $A_{w'(\ell)}^{2}$ the node factor
obtained from $A_{w'(\ell)}$ by replacing
the factor $({\rm e}^{-2\pi\ii\om\nu_{\ell}}-1)$
with $-2\pi\ii \, {\rm e}^{-2\pi\ii\om\nu_{\ell}}$.
The derivative of $\Val_{T}(u)$ with respect to $u=\om\nu_{\ell_{T}'}$
produce several contributions: for any line $\ell \in \calP_{T}$ there are
three terms which are obtained from $\Val_{T}(u)$, respectively,
(1) by replacing $\calG_{\ell}$ with $\partial_{u}\calG_{\ell}$,
(2) by replacing $A_{w(\ell)}$ with $A_{w(\ell)}^{1}$,
(3) by replacing $A_{w'(\ell)}$ with $A_{w'(\ell)}^{2}$,
To distinguish the three contributions we can introduce a label $\gamma=1,2,3$,
so that all the contributions produced by differentiation will be indexed
by two labels $(\ell,\gamma)$, with $\ell\in\calP_{T}$ and $\gamma\in\{1,2,3\}$
and denoted by $\calU_{T}(\om\nu;\ell,\gamma)$.

When considering $\partial_{u}^{2}\Val_{T}(u)$ in (\ref{eq:4.4}),
one has a sum of contributions which can be identified by four labels
$(\ell_{1},\ell_{2},\gamma_{1},\gamma_{2})$,
with $\ell_{1},\ell_{2}\in\calP_{T}$ and $\gamma_{1},\gamma_{2}\in\{1,2,3\}$,
and denoted by $\calU_{T}(u;\ell_{1},\gamma,\ell_{2},\gamma_{2})$;
for $\ell_{1}=\ell_{2}$ we shall impose the constraint $\gamma_{1} \le \gamma_{2}$ to avoid overcountings.
More precisely $\calU_{T}(u;\ell_{1},\gamma_1,\ell_{2},\gamma_1)$ is obtained from $\Val_{T}(u)$
by the following replacements
\begin{equation} \nonumber
\begin{array}{llll}
\ell_{1}=\ell_{2}, & \gamma_{1} = \gamma_{2} = 1, & &
\calG_{\ell_{1}} \to \partial_{u}^{2} \calG_{\ell_{1}} \\
\ell_{1}=\ell_{2}, & \gamma_{1} = 1, \gamma_{2} = 2 , & &
\calG_{\ell_{1}} A_{w(\ell_{1})} \to \partial_{u} \calG_{\ell_{1}} A^{1}_{w(\ell_{1})} \\
\ell_{1}=\ell_{2}, & \gamma_{1} = 1, \gamma_{2} = 3 , & &
\calG_{\ell_{1}} A_{w'(\ell_{1})} \to \partial_{u} \calG_{\ell_{1}} A^{2}_{w'(\ell_{1})} \\
\ell_{1}=\ell_{2}, & \gamma_{1} = \gamma_{2} = 2 , & &
A_{w(\ell_{1})} \to A^{3}_{w(\ell_{1})} \\
\ell_{1}=\ell_{2}, & \gamma_{1} = 2, \gamma_{2} = 3 , & &
A_{w(\ell_{1})} A_{w'(\ell_{1})} 
\to A^{1}_{w(\ell_{1})} A^{2}_{w'(\ell_{1})} \\
\ell_{1}=\ell_{2}, & \gamma_{1} = \gamma_{2} = 3 , & &
A_{w'(\ell_{1})} \to A^{4}_{w'(\ell_{1})} \\
\ell_{1} \neq \ell_{2}, & \gamma_{1} = \gamma_{2} = 1, & &
\calG_{\ell_{1}} \calG_{\ell_{2}} \to
\partial_{u} \calG_{\ell_{1}} \partial_{u} \calG_{\ell_{2}} \\
\ell_{1}\neq\ell_{2}, & \gamma_{1} = 1, \gamma_{2} = 2 , & &
\calG_{\ell_{1}} A_{w(\ell_{2})} \to \partial_{u} \calG_{\ell_{1}} A^{1}_{w(\ell_{2})} \\
\ell_{1}\neq\ell_{2}, & \gamma_{1} = 2, \gamma_{2} = 1 , & &
\calG_{\ell_{2}} A_{w(\ell_{1})} \to \partial_{u} \calG_{\ell_{2}} A^{1}_{w(\ell_{1})} \\
\ell_{1}\neq\ell_{2}, & \gamma_{1} = 1, \gamma_{2} = 3 , & &
\calG_{\ell_{1}} A_{w'(\ell_{2})} \to \partial_{u} \calG_{\ell_{1}} A^{2}_{w'(\ell_{2})} \\
\ell_{1}\neq\ell_{2}, & \gamma_{1} = 3, \gamma_{2} = 1 , & &
\calG_{\ell_{2}} A_{w'(\ell_{1})} \to \partial_{u} \calG_{\ell_{2}} A^{2}_{w'(\ell_{1})} \\
\ell_{1}\neq\ell_{2}, & \gamma_{1} = \gamma_{2} = 2 , & &
A_{w(\ell_{1})} A_{w(\ell_{2})} \to A^{1}_{w(\ell_{1})} A^{1}_{w(\ell_{2})} \\
\ell_{1}\neq\ell_{2}, & \gamma_{1} = 2, \gamma_{2} = 3 , & &
A_{w(\ell_{1})} A_{w'(\ell_{2})} 
\to A^{1}_{w(\ell_{1})} A^{2}_{w'(\ell_{2})} \\
\ell_{1}\neq\ell_{2}, & \gamma_{1} = 3, \gamma_{2} = 2 , & &
A_{w(\ell_{2})} A_{w'(\ell_{1})} 
\to A^{1}_{w(\ell_{2})} A^{2}_{w'(\ell_{1})} \\
\ell_{1}\neq \ell_{2}, & \gamma_{1} = \gamma_{2} = 3 , & &
A_{w'(\ell_{1})} A_{w'(\ell_{2})} \to A^{2}_{w'(\ell_{1})} A^{2}_{w'(\ell_{2})}
\end{array}
\end{equation}
where the node factor $A^{3}_{w(\ell)}$ is obtained from $A_{w(\ell)}$ 
by replacing the factor $\de_{+}(\om\nu_{\ell})$
with $(2\pi\ii)^{2} {\rm e}^{2\pi\ii\om\nu_{\ell}}$ if $\be_{\ell}=-$
and with $0$ if $\be_{\ell}=+$, and the node factor
$A^{4}_{w'(\ell)}$ is obtained from $A_{w'(\ell)}$ by replacing
the factor $({\rm e}^{-2\pi\ii\om\nu_{\ell}}-1)$
with $(-2\pi\ii)^{2} {\rm e}^{-2\pi\ii\om\nu_{\ell}}$.
The list of values above may look a bit tricky: however, the labels
$\ell_{1},\ell_{2},\gamma_1,\gamma_2$ simply identify the quantities
that have to be differentiated in $\Val_{T}(u)$. Therefore we have
\begin{subequations} \label{eq:4.6}
\begin{align}
\partial_{u} \Val_{T}(u) & = \sum_{\ell_{1}\in\calP_{T}} \sum_{\gamma_{1}=1,2,3}
\calU_{T}(u;\ell_{1},\gamma_{1}) ,
\label{eq:4.6a} \\
\partial_{u}^{2} \Val_{T}(u) & = \sum_{\ell_{1},\ell_{2}\in\calP_{T}}
\sideset{}{'}\sum_{\gamma_{1},\gamma_{2}=1,2,3}
\calU_{T}(u;\ell_{1},\gamma_{1},\ell_{2},\gamma_{2}) .
\label{eq:4.6b}
\end{align}
\end{subequations}
where, here and henceforth, the prime in the last sum recalls the constraint
$\gamma_{1}\le \gamma_{2}$ when $\ell_{1}=\ell_{2}$.
In (\ref{eq:4.6a}) we write
\begin{equation} \label{eq:4.7}
\calU_{T}(u;\ell_{1},\gamma_{1}) =
\Big( \prod_{v\in V(T)} \overline{A}_{v} \Big)
\Big( \prod_{\ell\in L(T)} \overline{\calG}_{\ell} \Big) ,
\end{equation}
where
\begin{equation} \nonumber
\overline{\calG}_{\ell} : =
\begin{cases}
\partial_{u}\calG_{\ell} , & \ell=\ell_{1}, \gamma_{1}=1 , \\
\calG_{\ell} , & \hbox{otherwise} ,
\end{cases}
\qquad
\overline{A}_{v}: = \begin{cases}
A_{v}^{1}, & v=w(\ell) , \gamma=2 , \\
A_{v}^{2} , & v=w'(\ell), \gamma=3 , \\
A_{v} , & \hbox{otherwise} .
\end{cases}
\end{equation}
Analogously we can write in (\ref{eq:4.6b})
\begin{equation} \label{eq:4.8}
\calU_{T}(t_{T}u;\ell_{1},\gamma_{1},\ell_{2},\gamma_{2}) =
\Big( \prod_{v\in V(T)} \overline{A}_{v} \Big)
\Big( \prod_{\ell\in L(T)} \overline{\calG}_{\ell} \Big) ,
\end{equation}
where
\begin{equation} \label{eq:4.9}
\overline{\calG}_{\ell} : = \begin{cases}
\partial_{u}^{2}\calG_{\ell}, & \ell=\ell_{1}=\ell_{2}, \gamma_{1}=\gamma_{2}=1 , \\
\partial_{u}\calG_{\ell} , & \ell=\ell_{1}, \gamma_{1}=1 \hbox{ or }
\ell=\ell_{2}, \gamma_{2}=1, \hbox{ with } \ell_{1}\neq\ell_{2}, \\
\calG_{\ell} , & \hbox{otherwise} ,
\end{cases}
\end{equation}
and
\begin{equation} \label{eq:4.10}
\overline{A}_{v} := \begin{cases}
A_{v}^{3} , &  v=w(\ell_{1})=w(\ell_{2}) , \gamma_{1}=\gamma_{2}=2 , \\
A_{v}^{4} , &  v=w'(\ell_{1})=w'(\ell_{2}), \gamma_{1}=\gamma_{2}=3 , \\
A_{v}^{1} , &  v=w(\ell_{1}), \gamma_{1}=2 \hbox{ or }
v=w(\ell_{2}), \gamma_{1}=2 , \hbox{ with } \ell_{1}\neq\ell_{2}, \\
A_{v}^{2} , &  v=w'(\ell_{1}), \gamma_{1}=3 \hbox{ or }
v=w'(\ell_{2}), \gamma_{1}=3 , \hbox{ with } \ell_{1}\neq\ell_{2}, \\
A_{v}  , & \hbox{otherwise} .
\end{cases}
\end{equation}

%%%%%%%%%%%%%%%%%%%%%%%%%%%%%%%%%%%%%%%%%%%%%%%%%
\begin{rmk} \label{rmk:4.2}
\emph{
The value of both $\overline{\calG}_{\ell}$ and $\overline{A}_{v}$ is well-defined
as it depends only on the number of derivatives acting on $\cal{G}_{\ell}$ and $A_{v}$,
and not to the fact that it appears in a contribution to the first or second
derivative of $\Val_{T}(u)$.
}
\end{rmk}
%%%%%%%%%%%%%%%%%%%%%%%%%%%%%%%%%%%%%%%%%%%%%%%%%

Given a self-energy cluster $T\in\gotT_{1}(\theta)$ with $\de_{T}=\LL$,
as in Section \ref{sec:3} we call $\gotF(T)$ the set of all self-energy clusters 
$T'\in\gotS_{k,n}$, where $k$ and $n$ are the order and scale, respectively, of $T$,
with $\nu_{\ell_{T'}'}=\nu_{\ell_{T}'}$.
The following result allows us to get rid of all contributions (\ref{eq:4.5})
in which at least one self-energy cluster $T\in\gotT_{1}(\theta)$ has $\de_{T}=\LL$.

%%%%%%%%%%%%%%%%%%%%%%%%%%%%%%%%%%%%%%%%%%%%%%%%%
\begin{lemma} \label{lem:4.3}
Let $T$ be a maximal self-energy cluster of $\theta$ of order $\overline{k}$
and scale $\overline{n}$. The sum of $\Val(\theta')$ over
all $\theta'\in\Theta_{k,\nu}$ obtained from $\theta$ by inserting
any $T'\in\gotF(T)$ instead of $T$ and replacing $\Val_{T'}(\om\nu_{\ell_{T'}'})$
with $\LL\Val_{T'}(\om\nu_{\ell_{T'}'})$ gives zero.
\end{lemma}
%%%%%%%%%%%%%%%%%%%%%%%%%%%%%%%%%%%%%%%%%%%%%%%%%

%%%%%%%%%%%%%%%%%%%%%%%%%%%%%%%%%%%%%%%%%%%%%%%%%
\prova The result follows from Lemma \ref{lem:3.16}.
\qed
%%%%%%%%%%%%%%%%%%%%%%%%%%%%%%%%%%%%%%%%%%%%%%%%%

\vspace{0.3cm}

So we are left only with the contribution (\ref{eq:4.5}) with all self-energy clusters
$T\in \gotT_{1}(\theta)$ having $\de_{T}=\RR$. 
We can write each factor $\RR \Val_{T}(\om\nu_{\ell_{T}'})$
in (\ref{eq:4.5}) according to (\ref{eq:4.4}), (\ref{eq:4.6b}) and (\ref{eq:4.8}), that is
\begin{equation} \label{eq:4.11}
\RR \Val_{T}(\om\nu_{\ell_{T}'}) = \sum_{\ell_{1},\ell_{2}\in\calP_{T}}
\sideset{}{'}\sum_{\gamma_{1},\gamma_{2}=1,2,3}
(\om\nu_{\ell_{T}'})^{2} \int_{0}^{1} {\rm d}t_{T} \left( 1 - t_{T} \right)
\Big( \prod_{v\in V(T)} \overline{A}_{v} \Big)
\Big( \prod_{\ell\in L(T)} \overline{\calG}_{\ell} \Big) .
\end{equation}

Given a line $\ell\in T$, define the \emph{cloud} of $\ell$ in $T$ as the set
$\gotC_{\ell}(T)= \{T_{1},\ldots,T_{p}\}$ such that
$T_{p}$ is the minimal self-energy cluster containing $\ell$ and
$T_{1} \supset T_{2} \supset \ldots \supset T_{p}$,
with $T_{j}\in\TT_{1}(T_{j-1})$ for $j=2,\ldots,p$ and $T_{1}\in\gotT_{1}(T)$.
For fixed $\ell_{1},\ell_{2}\in\calP_{T}$ let $\gotC_{\ell_{1}}(T)$
and $\gotC_{\ell_{2}}(T)$ the clouds of the two lines in $T$
(if the two lines coincide there is only one cloud).
We associate a label $\zeta_{T'}=0$ with each self-energy cluster
$T' \in \gotC_{\ell_{1}}(T) \cap \gotC_{\ell_{2}}(T)$.
We denote by $\gotT_{0}(T)$ the set of
such self-energy clusters and set $\gotT^{*}(T)=\gotT(T)\setminus\gotT_{0}(T)$.
Call $\gotT^{*}_{1}(T)$ the set of maximal self-energy clusters in $\gotT^{*}(T)$.
We associate with each $T'\in\gotT^{*}_{1}(T)$ a label $\zeta_{T'}=1$
if $T'$ contains one of the two lines $\ell_{1}$ and $\ell_{2}$
(by construction it cannot contain both of them) and a label $\zeta_{T'}=2$
if it does not contain any of them.
Finally denote by $T^{*}$ the set of nodes and lines in $T$
which are outside the self-energy clusters $T'\in\gotT^{*}_{1}(T)$.

%%%%%%%%%%%%%%%%%%%%%%%%%%%%%%%%%%%%%%%%%%%%%%%%%
\begin{rmk} \label{rmk:4.4}
\emph{
Note that both $\gotT_{1}^{*}(T)$ and $T^{*}$ depend on $\ell_{1}$ and $\ell_{2}$.
One can think of $\gotT_{1}^{*}(T)$ as the set of self-energy clusters which
become maximal in $T$ when ignoring the self-energy clusters
which belong to both the clouds of $\ell_{1}$ and $\ell_{2}$.
}
\end{rmk}
%%%%%%%%%%%%%%%%%%%%%%%%%%%%%%%%%%%%%%%%%%%%%%%%%

%%%%%%%%%%%%%%%%%%%%%%%%%%%%%%%%%%%%%%%%%%%%%%%%%
\begin{rmk} \label{rmk:4.5}
\emph{
If $T'\in \gotT^{*}_{1}(T)$ and $\zeta_{T'}=2$, one has
$\overline{\calG}_{\ell}=\calG_{\ell}$ for all lines $\ell\in L(T')$.
If $T'\in \gotT^{*}_{1}(T)$ and $\zeta_{T'}=1$, there is one line
$\ell_{1}\in L(T')$ such that $\overline{\calG}_{\ell_{1}}=\partial_{u}\calG_{\ell_{1}}$,
while $\overline{\calG}_{\ell}=\calG_{\ell}$ for $\ell\in L(T')\setminus\{\ell_{1}\}$.
}
\end{rmk}
%%%%%%%%%%%%%%%%%%%%%%%%%%%%%%%%%%%%%%%%%%%%%%%%%

By defining
\begin{subequations} \label{4.12}
\begin{align}
\overline{\Val}_{T^{*}}(\om\nu_{\ell_{T}'}(t_T))  & =
\Big( \prod_{v\in V(T^{*})} \overline{A}_{v} \Big)
\Big( \prod_{\ell\in L(T^{*})} \overline{\calG}_{\ell} \Big) ,
\label{eq:4.12a} \\
\overline{\Val}_{T'}(\om\nu_{\ell_{T'}'}(t_T))  & =
\Big( \prod_{v\in V(T')} \overline{A}_{v} \Big)
\Big( \prod_{\ell\in L(T')} \overline{\calG}_{\ell} \Big) , \qquad
T' \in \gotT^{*}_{1}(T) ,
\label{eq:4.12b}
\end{align}
\end{subequations}
where $\nu_{\ell}(t_{T})=\nu_{\ell}$ if $\ell\notin\calP_{T}$ and
$\nu_{\ell}(t_{T})=\nu_{\ell}+t_{T}\nu_{\ell_{T}'}$
if $\ell\in\calP_{T}$, we can rewrite (\ref{eq:4.11}) as
\begin{eqnarray} \label{eq:4.13}
& & \RR \Val_{T}(\om\nu_{\ell_{T}'}) = \sum_{\ell_{1},\ell_{2}\in\calP_{T}}
\sideset{}{'}\sum_{\gamma_{1},\gamma_{2}=1,2,3}
(\om\nu_{\ell_{T}'})^{2} \times \\
& & \qquad \qquad \times \int_{0}^{1} {\rm d}t_{T} \left( 1 - t_{T} \right)
\overline{\Val}_{T^{*}}(\om\nu_{\ell_{T}'}(t_T)) \prod_{T'\in\gotT^{*}_{1}(T)}
\overline{\Val}_{T'}(\om\nu_{\ell_{T'}'}(t_T))  . \nonumber
\end{eqnarray}

For $T' \in \gotT^{*}_{1}(T)$ we define
\begin{equation} \label{eq:4.14}
\LL \overline{\Val}_{T'}(u) = \begin{cases}
\overline{\Val}_{T'}(0) + u \, \partial_{u} \overline{\Val}_{T'}(0) , & \zeta_{T'}=2 , \\
\overline{\Val}_{T'}(0) , & \zeta_{T'}=1 , 
\end{cases}
\end{equation}
and set $\RR \overline{\Val}_{T'}(u) = \overline{\Val}_{T'}(u) - \LL \overline{\Val}_{T'}(u)$, so that
\begin{equation} \label{eq:4.15}
\begin{cases}
\displaystyle{ \RR\overline{\Val}_{T'}(u) = u^{2} \int_{0}^{1} {\rm d}t_{T'} \left( 1 - t_{T'} \right)
\partial_{u}^{2} \overline{\Val}_{T'}(t_{T'}u) } , & \qquad \zeta_{T'} =2 , \\

\displaystyle{ \RR\overline{\Val}_{T'}(u) = u \int_{0}^{1} {\rm d}t_{T'} \,
\partial_{u} \overline{\Val}_{T'}(tu) } , & \qquad \zeta_{T'} =1 .
\end{cases}
\end{equation}

Then we repeat the construction by associating a label $\de_{T'}\in\{\LL,\RR\}$
with each $T'\in\gotT^{*}_{1}(T)$ and write $\RR \Val_{T}(\om\nu_{\ell_{T}'})$
in (\ref{eq:4.13}) as sum of contributions
\begin{eqnarray} \label{eq:4.16}
& & \sum_{\ell_{1},\ell_{2}\in\calP_{T}}
\sideset{}{'}\sum_{\gamma_{1},\gamma_{2}=1,2,3}
(\om\nu_{\ell_{T}'})^{2} \int_{0}^{1} {\rm d}t_{T} \left( 1 - t_{T} \right) \times \\
& & \qquad \qquad \times
\overline{\Val}_{T^{*}}(\om\nu_{\ell_{T}'}(t_T))
\Big( \prod_{\substack{ T'\in\gotT^{*}_{1}(T) \\ \de_{T}=\RR}}
\RR\overline{\Val}_{T'}(\om\nu_{\ell_{T'}'}(t_T)) \Big)
\Big( \prod_{\substack{ T'\in\gotT^{*}_{1}(T) \\ \de_{T}=\LL}}
\LL\overline{\Val}_{T'}(\om\nu_{\ell_{T'}'}(t_T)) \Big) ,
\nonumber
\end{eqnarray}
where each $\RR\overline{\Val}_{T'}(\om\nu_{\ell_{T'}'}(t_T))$
is expressed as in (\ref{eq:4.15}), with $\partial_{u}^{2} \overline{\Val}_{T'}(t_{T'}u)$
and $\partial_{u} \overline{\Val}_{T'}(t_{T'}u)$
written as in (\ref{eq:4.6}) to (\ref{eq:4.8}). In particular the momentum of
any line $\ell\in L(T')$ becomes
\begin{equation} \label{eq:4.17}
\nu_{\ell}(t_{T'},t_{T}) := \begin{cases}
\nu_{\ell}^{0} , & 
\ell\notin\calP_{T'} , \\
\nu_{\ell}^{0}+t_{T'}\nu_{\ell_{T'}'}^{0} , & 
\ell\in\calP_{T'} \hbox{ and } \ell_{T'}'\notin\calP_{T} , \\
\nu_{\ell}^{0}+t_{T'}\,(\nu_{\ell_{T'}'}^{0}+t_{T}\nu_{\ell_{T}}) , &
\ell\in\calP_{T'} \hbox{ and }\ell_{T'}'\notin\calP_{T} .
\end{cases}
\end{equation}
The contributions (\ref{eq:4.16}) in which at least one self-energy cluster
$T'\in\gotT^{*}_{1}(T)$ has $\de_{T'}=\LL$ vanish when summed together.
This follows from the following result.

%%%%%%%%%%%%%%%%%%%%%%%%%%%%%%%%%%%%%%%%%%%%%%%%%
\begin{lemma} \label{lem:4.6}
Let $\overline{k}$ and $\overline{n}$ be the order and the scale, respectively,
of $T'\in\gotT_{1}^{*}(T)$. The sum of $\Val(\theta')$ over all $\theta'\in
\Theta_{k,\nu}$ obtained from $\theta$ by inserting any self-energy cluster
$T''\in\gotF(T')$ instead of $T'$ and replacing
$\Val_{T''}(\om\nu_{\ell_{T''}'})$ with $\LL\Val_{T''}(\om\nu_{\ell_{T''}'})$ gives zero.
\end{lemma}
%%%%%%%%%%%%%%%%%%%%%%%%%%%%%%%%%%%%%%%%%%%%%%%%%

%%%%%%%%%%%%%%%%%%%%%%%%%%%%%%%%%%%%%%%%%%%%%%%%%
\prova One reasons as for Lemma \ref{lem:4.3}.
\qed
%%%%%%%%%%%%%%%%%%%%%%%%%%%%%%%%%%%%%%%%%%%%%%%%%

\vspace{0.3cm}

Therefore we have
\begin{eqnarray} \label{eq:4.18}
& & \sum_{\theta\in\Theta_{k,n}} \Val(\theta) =
\sum_{\theta\in\Theta_{k,n}} \Val(\mathring{\theta})
\prod_{T\in\gotT_{1}(\theta)} \sum_{\ell_{T,1},\ell_{T,2}\in\calP_{T}} 
\sideset{}{'}\sum_{\gamma_{T,1},\gamma_{T,2}=1,2,3}
(\om\nu_{\ell_{T}'})^{2} \times \nonumber \\
& & \qquad \qquad \times \int_{0}^{1} {\rm d}t_{T} \left( 1 - t_{T} \right)
\overline{\Val}_{T^{*}}(\om\nu_{\ell_{T}'}(t_T)) \prod_{T'\in\gotT^{*}_{1}(T)}
\RR \overline{\Val}_{T'}(\om\nu_{\ell_{T'}'}(t_T)) .
\end{eqnarray}
The construction can be iterated further.
At each step the order of the self-energy clusters has
decreased, so that eventually the procedure stops.
Moreover, every time we split the value of the a self-energy cluster
into the sum of the localised value plus the regularised value,
we can neglect the localised value. Indeed the following extension
of Lemma \ref{lem:4.6} holds.

%%%%%%%%%%%%%%%%%%%%%%%%%%%%%%%%%%%%%%%%%%%%%%%%%
\begin{lemma} \label{lem:4.7}
Given any $T\in\gotT(T)$ denote by $\overline{k}$ and $\overline{n}$
the order and the scale, respectively, of $T$.
The sum of $\Val(\theta')$ over all $\theta'\in\Theta_{k,\nu}$
obtained from $\theta$  by inserting
any $T'\in\gotF(T)$ instead of $T$ and replacing $\Val_{T'}(\om\nu_{\ell_{T'}'})$
with $\LL\Val_{T'}(\om\nu_{\ell_{T'}'})$ gives zero.
\end{lemma}
%%%%%%%%%%%%%%%%%%%%%%%%%%%%%%%%%%%%%%%%%%%%%%%%%

%%%%%%%%%%%%%%%%%%%%%%%%%%%%%%%%%%%%%%%%%%%%%%%%%
\prova Simply note that, when computing the localised value of a self-energy cluster $T$,
one puts $\nu_{\ell}=\nu_{\ell}^{0}$ for all lines $\ell\in L(T)$.
Therefore the cancellation mechanisms
underlying Lemma \ref{lem:3.16} apply to any self-energy cluster,
indepenedently of the fact that it is maximal or inside other self-energy clusters.
\qed
%%%%%%%%%%%%%%%%%%%%%%%%%%%%%%%%%%%%%%%%%%%%%%%%%

\vspace{0.3cm}

With any self-energy cluster $T$ we associate a label $\zeta_{T}\in\{0,1,2\}$,
an interpolation parameter $t_{T}\in[0,1]$ and a measure 
$\pi_{\zeta_{T}}(t_{T})\,{\rm d}t_{T}$, where
\begin{equation} \label{eq:4.19}
\pi_{\zeta}(t) = \begin{cases}
\left( 1 - t \right) , & \zeta = 2, \\
1 , & \zeta=1 , \\
\delta(t-1) , & \zeta=0 ,
\end{cases}
\end{equation}
where $\delta$ is the Dirac delta. 

Denote by $\bt=\{t_{T}\}$ the set of all
interpolation parameters and define recursively
\begin{equation} \label{eq:4.20}
\nu_{\ell}(\bt) =
\begin{cases}
\nu_{\ell}^{0} , & \ell \notin \calP_{T} , \\
\nu_{\ell}^{0} + t_{T} \nu_{\ell_{T}'}(\bt) , & \ell \in \calP_{T} ,
\end{cases}
\end{equation}
where $T$ is the minimal self-energy cluster containing $\ell$ (if any).

%%%%%%%%%%%%%%%%%%%%%%%%%%%%%%%%%%%%%%%%%%%%%%%%%
\begin{rmk} \label{rmk:4.8}
\emph{
For any $\ell\in L(\theta)$ the momentum $\nu_{\ell}(\bt)$ depends only
on the interpolation parameters of the self-energy clusters $T'$ such that
$\ell\in\calP_{T'}$.
}
\end{rmk}
%%%%%%%%%%%%%%%%%%%%%%%%%%%%%%%%%%%%%%%%%%%%%%%%%

Eventually one obtains that
\begin{eqnarray} \label{eq:4.21}
& & \sum_{\theta\in\Theta_{k,\nu}} \Val(\theta) =  
\sum_{\theta\in\Theta_{k,\nu}} 
\Big( \prod_{v\in V(\mathring{\theta})} \overline{A}_{v} \Big)
\Big( \prod_{\ell \in L(\mathring{\theta})} \overline{\calG}_{\ell} \Big) \times \\ 
& & \qquad \qquad \times \prod_{T\in\gotT(\theta)} {\sum}_{T}^{\zeta_{T}}
(\om\nu_{\ell_{T}'}(\bt))^{\zeta_{T}} \int_{0}^{1} \pi_{\zeta_{T}}(t_{T})\,{\rm d}t_{T}
\Big( \prod_{v\in V(T^{*})} \overline{A}_{v} \Big)
\Big( \prod_{\ell \in L(T^{*})} \overline{\calG}_{\ell} \Big)
\nonumber
\end{eqnarray}
where $\overline{\calG}_{\ell}$ and $\overline{A}_{v}$ are defined in (\ref{eq:4.9})
and (\ref{eq:4.10}), respectively, and we have introduced the shorthand notation
\begin{equation} \nonumber
{\sum}_{T}^{\zeta_{T}} =
\begin{cases}
\displaystyle{
\sum_{\ell_{T,1},\ell_{T,2}\in\calP_{T}} 
\sideset{}{'}\sum_{\gamma_{T,1},\gamma_{T,2}=1,2,3}} ,
& \zeta_{T}=2 , \\
\displaystyle{ \sum_{\ell_{T,1} \in\calP_{T}} \sum_{\gamma_{T,1}=1,2,3}} ,
& \zeta_{T}=1 .
\end{cases}
\end{equation}

%%%%%%%%%%%%%%%%%%%%%%%%%%%%%%%%%%%%%%%%%%%%%%%%%
\begin{rmk} \label{rmk:4.9}
\emph{
The propagator of each line is differentiated at most twice.
}
\end{rmk}
%%%%%%%%%%%%%%%%%%%%%%%%%%%%%%%%%%%%%%%%%%%%%%%%%

In (\ref{eq:4.21}) we can bound, for any $\xi_{2}\in(0,\xi)$,
\begin{equation} \label{eq:4.22}
\Big( \prod_{v\in V(\mathring{\theta})} \left| \overline{A}_{v} \right| \Big)
\prod_{T\in \gotT(\theta)}
\Big( \prod_{v\in V(T^{*})} \left| \overline{A}_{v} \right| \Big) 
\le \left( C_0 \Xi \rho^{-2} \xi_{2}^{-2} \right)^{k} 
\prod_{v\in V(\theta)} {\rm e}^{-(\xi-\xi_{2})|\nu_{v}|} ,
\end{equation}
and, if we denote by $L^{\bullet}(\theta)$ the set of non-resonant lines in $\theta$,
\begin{eqnarray} \label{eq:4.23}
& & \hskip-.5truecm
\Big( \prod_{\ell \in L(\mathring{\theta})} \left| \overline{\calG}_{\ell} \right| \Big)
\prod_{T \in \gotT(\theta)}
\Big( \prod_{\ell \in L(T^{*})} \left| \overline{\calG}_{\ell} \right| \Big) \nonumber \\
& & \le C_{1}^{k}
\Big( \prod_{\ell \in L^{\bullet}(\theta)} \|\om\nu_{\ell}(\bt)\|^{-2} \Big)
\Big( \prod_{T \in \gotT(\theta)} \|\om\nu_{\ell_{T}}(\bt)\|^{-2} \Big) \times \\
& & \times \Big( \prod_{\substack{T \in \gotT(\theta) \\ \zeta_{T}=2}} \|\om\nu_{\ell_{T}'}(\bt)\|^{2}
\|\om\nu_{\ell_{T,1}}(\bt)\|^{-1} \|\om\nu_{\ell_{T,2}}(\bt)\|^{-1} \Big)
\Big( \prod_{\substack{T \in \gotT(\theta) \\ \zeta_{T}=1}} \|\om\nu_{\ell_{T}'}(\bt)\|
\|\om\nu_{\ell_{T,1}}(\bt)\|^{-1} \Big) \nonumber
\end{eqnarray}
for suitable constants $C_0$ and $C_{1}$. 
For each line $\ell_{T,i}$ in (\ref{eq:4.23}) we can write
\begin{equation} \label{eq:4.24}
\|\om\nu_{\ell_{T,1}}(\bt)\|^{-1} = 
\|\om\nu_{\ell_{T,1}}(\bt)\|^{-1} \prod_{T \in \gotC_{\ell}(T_{\ell})}
\|\om\nu_{\ell_{T}'}(\bt)\|^{-1} \|\om\nu_{\ell_{T}'}(\bt)\| ,
\end{equation}
where $T_{\ell}$ the minimal self-energy cluster $T'$
containing $\ell_{T,i}$ with $\zeta_{T'} \neq 0$.

%%%%%%%%%%%%%%%%%%%%%%%%%%%%%%%%%%%%%%%%%%%%%%%%%
\begin{rmk} \label{rmk:4.10}
\emph{
Consider any summand in the right hand side of (\ref{eq:4.21}).
For each line $\ell$ with $n_{\ell}\ge 1$ the propagator $\overline{\calG}_{\ell}$
can be bounded proportionally to $\|\om\nu_{\ell}(\bt)\|^{-2-p}$,
if $\overline{\calG}_{\ell}=\partial_{u}^{p}\calG_{\ell}$ for $p\in\{0,1,2\}$,
where $1/16q_{n_{\ell}+1} \le \|\om\nu_{\ell}(\bt)\| \le 1/64q_{n_{\ell}}$
by Lemma \ref{lem:3.2}.
}
\end{rmk}
%%%%%%%%%%%%%%%%%%%%%%%%%%%%%%%%%%%%%%%%%%%%%%%%%

\vspace{.3cm}

In the light of Remark \ref{rmk:4.10}, if we use (\ref{eq:4.24}) and the fact
that $n_{\ell} \le n_{T}$ for all $\ell\in L(T)$, we can bound in (\ref{eq:4.23})
\begin{eqnarray} \label{eq:4.25}
& & \Big( \prod_{\substack{T \in \gotT(\theta) \\ \zeta_{T}=2}}
\|\om\nu_{\ell_{T}'}(\bt)\|^{2}
\|\om\nu_{\ell_{T,1}}(\bt)\|^{-1} \|\om\nu_{\ell_{T,2}}(\bt)\|^{-1} \Big)
\Big( \prod_{\substack{T \in \gotT(\theta) \\ \zeta_{T}=1}}
\|\om\nu_{\ell_{T}'}(\bt)\|
\|\om\nu_{\ell_{T,1}}(\bt)\|^{-1} \Big) \nonumber \\
& & \qquad \qquad \le C_{2}^{k}
\Big( \prod_{T\in\gotT(\theta)} \|\om\nu_{\ell_{T}'}(\bt)\|^{2} \Big)
\Big( \prod_{T\in\gotT(\theta)} q_{n_{T}+1}^{2} \Big) ,
\end{eqnarray}
for a suitable constant $C_{2}$, so as to obtain
\begin{equation} \label{eq:4.26}
\Big( \prod_{\ell \in L(\mathring{\theta})} \left| \overline{\calG}_{\ell} \right| \Big)
\prod_{T \in \gotT(\theta)}
\Big( \prod_{\ell \in L(T^{*})} \left| \overline{\calG}_{\ell} \right| \Big)
\le C_{1}^{k} C_{2}^{k}  \Big( \prod_{\ell \in L^{\bullet}(\theta)}
\|\om\nu_{\ell}(\bt)\|^{-2} \Big)
\Big( \prod_{T\in\gotT(\theta)} q_{n_{T}+1}^{2} \Big) .
\end{equation}
Therefore one has
\begin{equation} \label{eq:4.27}
\Big( \prod_{\ell \in L(\mathring{\theta})} \left| \overline{\calG}_{\ell} \right| \Big)
\prod_{T \in \gotT(\theta)}
\Big( \prod_{\ell \in L(T^{*})} \left| \overline{\calG}_{\ell} \right| \Big)
\le C_{3}^{k} \prod_{n\ge 0} q_{n+1}^{2(N_{n}^{\bullet}(\theta)+P_{n}(\theta))} .
\end{equation}

To bound $N^{\bullet}(\theta)+P_{n}(\theta)$ in (\ref{eq:4.27}) we shall use
Lemma \ref{lem:3.13}. However, we have first to check that the lines satisfy
the support properties. If the momentum of each line $\ell$ were $\nu_{\ell}$
this would follow immediately from Lemma \ref{lem:3.2}
(see Remark \ref{rmk:3.7}). Though,
the momenta are $\nu_{\ell}(\bt)$, and the compact support functions
$\Psi_{n_{\ell}}$ only assure that $\Psi_{n_{\ell}}(\|\om\nu_{\ell}(\bt)\|) \neq0$,
that is $1/64q_{n_{\ell}+1} \le \|\om\nu_{\ell}(\bt)\| \le 1/16q_{n_{\ell}}$.
We want to show that the latter property implies the following result.

%%%%%%%%%%%%%%%%%%%%%%%%%%%%%%%%%%%%%%%%%%%%%%%%%
\begin{lemma} \label{lem:4.11}
In each summand on the right hand side of (\ref{eq:4.21}), for all lines $\ell \in L(\theta)$
one has $1/128q_{n_{\ell}+1} \le \|\om\nu_{\ell}\| \le 1/8q_{n_{\ell}}$.
\end{lemma}
%%%%%%%%%%%%%%%%%%%%%%%%%%%%%%%%%%%%%%%%%%%%%%%%%

%%%%%%%%%%%%%%%%%%%%%%%%%%%%%%%%%%%%%%%%%%%%%%%%%
\prova We say that a line $\ell$ has depth $0$ if it is outside any self-energy cluster
and depth $p\ge 1$ if there are $p$ self-energy clusters
$T_{1} \supset T_{2} \supset \ldots \supset T_{p}$,
such that $T_{i} \in \gotT_{1}^{*}(T_{i-1})$ for $i=2,\ldots,p$,
$T_{1}\in\gotT_{1}(\theta)$ and $\ell\in T_{p}^{*}$.
In such a case one has $\nu_{\ell}(\bt)=\nu_{\ell}^{0}$ if $\ell\notin\calP_{T_{p}}$
and $\nu_{\ell}(\bt)=\nu_{\ell}^{0}+t_{T_{p}}\nu_{\ell_{T_{p}}'}(\bt)$ if
$\ell\in\calP_{T_{p}}$.

We want to prove the bound by induction on the depth of the lines.

For any line $\ell$ of depth $0$ one has $\nu_{\ell}(\bt)=\nu_{\ell}$, so that
the bound trivially holds.

Let $\ell$ be a line with depth $1$ and let $T$ be the self-energy cluster
containing $\ell$. Again if $\nu_{\ell}^{0}=\nu_{\ell}$ there is nothing to prove,
so we have to consider explicitly only the case
$\nu_{\ell}(\bt)=\nu_{\ell}^{0}+t_{T}\om\nu_{\ell_{T}'}$.
Let $n$ be the scale of $\ell_{T}'$. By definition of self-energy cluster one has
$|\nu_{\ell}^{0}| \le M(T)<q_{n}$ and hence $\|\om\nu_{\ell}^{0}\|>1/4q_{n}$.
On the other hand one has $\|\om\nu_{\ell_{T}'}\|\le 1/16q_{n}$, so that
\begin{equation} \label{eq:4.28}
\frac{1}{2} \left\| \om\nu_{\ell}^{0} \right\| \le
\|\om\nu_{\ell}^{0}\| - \|\om\nu_{\ell_{T}'}\| \le
\| \om\nu_{\ell} \| \le
\|\om\nu_{\ell}^{0}\| + \|\om\nu_{\ell_{T}'}\| \le
\frac{3}{2} \| \om\nu_{\ell}^{0} \| .
\end{equation}
Let $n_{\ell}$ be the scale of the line $\ell$: then
$1/64q_{n_{\ell}+1} \le \|\om\nu_{\ell}(\bt)\|\le 1/16q_{n_{\ell}}$
because $\Psi_{n_{\ell}}(\|\om\nu_{\ell}(\bt)\|)\neq0$
(see also Lemma \ref{lem:3.2}). The quantity $\|\om\nu_{\ell}(\bt)\|$ is between
$\|\om\nu_{\ell}^{0}\|$ and $\|\om\nu_{\ell}\|$.
If  $\|\om\nu_{\ell}^{0}\|<\|\om\nu_{\ell}\|$ one has
$\|\om\nu_{\ell}\| \ge \|\om\nu_{\ell}(\bt)\| \ge 1/64q_{n_{\ell}+1} \ge 1/128q_{n_{\ell}+1}$
and $\|\om\nu_{\ell}\| \le 2\|\om\nu_{\ell}^{0}\| \le 2 \|\om\nu_{\ell}(\bt)\| \le 1/8q_{n_{\ell}}$;
if  $\|\om\nu_{\ell}^{0}\| \ge \|\om\nu_{\ell}\|$ one has
$\|\om\nu_{\ell}\| \le \|\om\nu_{\ell}(\bt)\| \le 1/16q_{n_{\ell}} \le 1/8q_{n_{\ell}}$
and $\|\om\nu_{\ell}\| \ge \|\om\nu_{\ell}^{0}\|/2 \ge \|\om\nu_{\ell}(\bt)\|/2 \ge 1/128q_{n_{\ell}}$.
Therefore in both cases one has
$1/8q_{n_{\ell}} \ge \|\om\nu_{\ell}\| \ge 1/128q_{n_{\ell}+1}$.

Now let us assume that the bound holds for all lines of depth $\le k$
and show that then it holds also for lines of depth $k+1$. Let $\ell$ be
one of such lines and let $T$ be the minimal self-energy cluster with $\zeta_{T}\neq0$
containing $\ell$. Once more if $\nu_{\ell}(\bt)=\nu_{\ell}^{0}$ the bound
holds trivially. If instead one has $\nu_{\ell}(\bt)=\nu_{\ell}^{0}+t_{T}\nu_{\ell_{T}'}(\bt)$
one can reason as in the previous case.
One has $\nu_{\ell}=\nu_{\ell}^{0}+\nu_{\ell_{T}'}$ and $\ell_{T}'$ has depth $k$,
so that, by the inductive hypothesis, the bound $1/8q_{n} \ge
\|\om\nu_{\ell_{T}'}\| \ge 1/128q_{n+1}$ holds, if $n$ is the scale
of $\ell_{T}'$. On the other hand $|\nu_{\ell}^{0}|\le M(T)<q_{n}$ and hence
$\|\om\nu_{\ell}^{0}\|>1/4q_{n}$, so that (\ref{eq:4.28}) holds also for the line $\ell$.
Moreover $\|\om\nu_{\ell}(\bt)\|$, for all values of the interpolation parameters $\bt$,
is between $\|\om\nu_{\ell}^{0}\|$ and $\|\om\nu_{\ell}\|$.
Therefore one can apply the previous argument and bound $\|\om\nu_{\ell}\|$
in terms of the scale $n_{\ell}$ of the line $\ell$ so as to obtain the desired bound.
\qed
%%%%%%%%%%%%%%%%%%%%%%%%%%%%%%%%%%%%%%%%%%%%%%%%%

\vspace{.3cm}

By combining Lemma \ref{lem:3.13} with Lemma \ref{lem:4.11} we obtain,
for arbitrary $n_{0}\in\NNN$,
\begin{equation} \nonumber
C_{3}^{k} \prod_{n\ge 0} q_{n+1}^{2(N_{n}^{\bullet}(\theta)+P_{n}(\theta))} 
\le C_{3}^{k} q_{n_{0}}^{2k}
\Big( \prod_{n \ge n_{0}} q_{n+1}^{4M(\theta)/q_{n}} \Big) \le
C_{3}^{k} q_{n_{0}}^{2k} \exp \Big( 4 M(\theta) \sum_{n\ge n_{0}}
\frac{1}{q_{n}} \log q_{n+1} \Big) .
\end{equation}
The sum converges by the assumption that $\gotB(\om)<\io$, so that
for any $\xi_{3}\in(0,\x-\xi_{2})$ one can choose $n_{0}=n_{0}(\xi_{3})$
suitably large so that
\begin{equation} \label{eq:4.30}
\Big( \prod_{\ell \in L(\mathring{\theta})} \left| \overline{\calG}_{\ell} \right| \Big)
\prod_{T \in \gotT(\theta)}
\Big( \prod_{\ell \in L(T^{*})} \left| \overline{\calG}_{\ell} \right| \Big) \le
C_{3}^{k} C_{4}^{k}(\xi_{3})\, {\rm e}^{\xi_{3} M(\theta)} ,
\end{equation}
with $C_{4}(\xi_{3})=q_{n_{0}(\xi_{3})}^{2}$.

By collecting together the bounds (\ref{eq:4.22}) and (\ref{eq:4.30}),
and taking $\xi_{1}\in(0,\xi-\xi_{2}-\xi_{3})$, we obtain
that each summand in the right hand side of (\ref{eq:4.21}) is bounded by
\begin{equation} \label{eq:4.31}
\left( \Xi \, C_{0} C_{3} \rho^{-2}\xi_{2}^{-2}  C_{4}(\xi_{3})\right)^{k}  {\rm e}^{-\xi_{1} |\nu|}
\prod_{v\in V(\theta)} {\rm e}^{-(\xi-\xi_{1}-\xi_{2}-\xi_{3}) |\nu_{v}|} .
\end{equation} 
We have still to perform the sum over the tree labels.
This is controlled by the following result.

%%%%%%%%%%%%%%%%%%%%%%%%%%%%%%%%%%%%%%%%%%%%%%%%%
\begin{lemma} \label{lem:4.12}
There exists a positive constant $C_{5}$ such that the
the sum of the quantities (\ref{eq:4.31}) over the tree labels is bounded by
$C_{5}^{k}C_{4}^{k}(\xi_{3})(\xi-\xi_{1}-\xi_{2}-\xi_{3})^{-k}\xi_{2}^{-2k}{\rm e}^{-\xi_{1}|\nu|}$.
\end{lemma}
%%%%%%%%%%%%%%%%%%%%%%%%%%%%%%%%%%%%%%%%%%%%%%%%%

%%%%%%%%%%%%%%%%%%%%%%%%%%%%%%%%%%%%%%%%%%%%%%%%%
\prova
The number of unlabelled trees of order $k$ is bounded by $4^{k}$. 
Then we have to sum over all the labels. The sum over the
labels $\zeta_{T}$ and the choices of the lines $\ell_{T,1},\ell_{T,2}$ if $\zeta_{T}=2$
and $\ell_{T,1}$ if $\zeta_{T}=1$ is bounded
by $C_{6}^{k}$ for some positive constant $C_{6}$. The sum over the Fourier labels
can be bounded thanks to the factors ${\rm e}^{-(\xi-\xi_{1}-\xi_{2})|\nu_{v}|}$
associated with the nodes and produces a further factor $C_{7}^{k}(\xi-\xi_{1}-\xi_{2}-\xi_{3})^{-k}$,
for a suitable constant $C_{7}$. Then the assertion follows
with $C_{5}=\Xi \rho^{-2} C_{0} C_{3} C_{6} C_{7}$.
\qed
%%%%%%%%%%%%%%%%%%%%%%%%%%%%%%%%%%%%%%%%%%%%%%%%%

\vspace{.3cm}

This concludes the proof of Theorem \ref{thm:2.0} when $\ols_{\e}(x,z)$ does not depend
explicitly on $\e$. If on the contrary $\ols_{\e}(x,z)$ depends analytically
also on $\e$, in (\ref{eq:2.2}) one has
\begin{equation} \nonumber
\ols_{\e}(x,z) = \sum_{n\in\ZZZ} \sum_{s=0}^{\io} {\rm e}^{\ii\n x} 
\s_{\n}(z,\e) , \qquad \frac{1}{s!q!} \left| \partial_{z}^{q} \partial_{\e}^{s}
\s_{\n}(2\pi\om,0) \right| \le \Xi {\rm e}^{-\xi |\n|} \rho^{-q} \mu^{-s} ,
\end{equation}
for suitable constants $\Xi$, $\xi$, $\rho$ and $\mu$,
provided $|\e|<\e_{1}$ for some $\e_{1}>0$. The discussion then proceeds as above,
the only difference being that each node carries a further label $s_{v}\in \ZZZ_{+}$ and
the node factor is given by (\ref{eq:2.11}), with $\s_{\n_v}(2\pi\om)$ replaced with
$(1/s_{v}!)\partial_{\e}^{s_{v}}\s_{\n_v}(2\pi\om,0)$.
The order of a tree $\theta$ is defined as
%
%\vspace{-.2cm}
\begin{equation} \nonumber
k= \sum_{v \in V(\theta)} \left( 1 + s_{v} \right) .
\end{equation}
Then, up to obvious minor changes, everything goes through
as before and the same results hold, simply with a different constant $\e_0$.

%%%%%%%%%%%%%%%%%%%%%%%%%%%%%%%%%%%%%%%%%%%%%%%%%
%%%%%%%%%%%%%%%%%%%%%%%%%%%%%%%%%%%%%%%%%%%%%%%%%
\zerarcounters 
\section{Some comments}
\label{sec:5} 
%%%%%%%%%%%%%%%%%%%%%%%%%%%%%%%%%%%%%%%%%%%%%%%%%
%%%%%%%%%%%%%%%%%%%%%%%%%%%%%%%%%%%%%%%%%%%%%%%%%

The results of Sections \ref{sec:4} can be summarised in the statement
of Theorem \ref{thm:2.0}. The result can be strengthened into Theorem \ref{thm:1.2},
as the following argument shows.

\vskip.4truecm

%%%%%%%%%%%%%%%%%%%%%%%%%%%%%%%%%%%%%%%%%%%%%%%%%
\noindent{\it Proof of Theorem \ref{thm:1.2}}.
First note that Theorem \ref{thm:2.0} ensures that if $\om$ satisfies the Bryuno condition
then for $\e$ small enough there exists an invariant curve with rotation number $\om$.
So, to complete the proof we have to show that if $\om$ does not satisfy the
Bryuno condition, then there exists a symplectic map $\Phi$ arbitrarily close
to $\Phi_0$ such that no analytic invariant curve with rotation number $\om$
and analytic in $\e$ exists. One can take as $\Phi$ the standard map.
\begin{equation} \label{eq:5.1}
\begin{cases}
x' = x + y + \e \, K \, \sin x , \\
y' = y + \e \, K \, \sin x ,
\end{cases}
\end{equation}
where a suitable constant $K$ has been introduced in order to make $|\Phi_{1}|=1$ in $\calD$.
Indeed, for such a map, if one denotes by $\rho_{0}(\om)$ the radius of analyticity
of the conjugation, one has \cite{D,BG}
\begin{equation} \label{eq:5.2}
\left| \log \rho_{0}(\om) + 2\gotB(\om) \right| \le C ,
\end{equation}
for some universal constant $C$. This means that the radius of convergence
vanishes when $\om$ does not satisfy the Bryuno condition: hence no
analytic invariant curve analytic in $\e$ exists in such a case.
\qed
%%%%%%%%%%%%%%%%%%%%%%%%%%%%%%%%%%%%%%%%%%%%%%%%%

%%%%%%%%%%%%%%%%%%%%%%%%%%%%%%%%%%%%%%%%%%%%%%%%%
\begin{rmk} \label{rmk:5.3}
\emph{
Of course, for the existence of an analytic invariant curve analytic in $\e$,
the Bryuno condition cannot be optimal for every perturbation.
A trivial counterexample is the null perturbation for which every invariant curve
exists independently of its rotation number. More generally one can consider
a perturbation which vanishes on a fixed unperturbed curve. Simply take $\Phi$
given by (\ref{eq:2.1}) with $\s_{\e}(x,x')=(x'-x-y_0)^{2}F_{\e}(x,x'-x)$,
with $F_{\e}$ an arbitrary analytic function $2\pi$-periodic in the first argument:
the invariant curve $x'=x+y_0$, $y'=y=y_0$ persists
for any $F_{\e}$, without any further assumption on the rotation number $\om=y_0/2\pi$.
}
\end{rmk}
%%%%%%%%%%%%%%%%%%%%%%%%%%%%%%%%%%%%%%%%%%%%%%%%%

%%%%%%%%%%%%%%%%%%%%%%%%%%%%%%%%%%%%%%%%%%%%%%%%%
\begin{rmk} \label{rmk:5.4}
\emph{
An interesting problem is whether the Bryuno condition is optimal for the existence 
of an analytic invariant curve analytic in $\e$ for generic perturbations of $\Phi_{0}$.
}
\end{rmk}
%%%%%%%%%%%%%%%%%%%%%%%%%%%%%%%%%%%%%%%%%%%%%%%%%

%%%%%%%%%%%%%%%%%%%%%%%%%%%%%%%%%%%%%%%%%%%%%%%%%
\begin{rmk} \label{rmk:5.5}
\emph{
Another interesting problem is whether a result like Theorem \ref{thm:1.2}
extends to the more general systems of the form (\ref{eq:1.1}),
which will be considered in Section \ref{sec:6}.
}
\end{rmk}
%%%%%%%%%%%%%%%%%%%%%%%%%%%%%%%%%%%%%%%%%%%%%%%%%

%%%%%%%%%%%%%%%%%%%%%%%%%%%%%%%%%%%%%%%%%%%%%%%%%
\begin{rmk} \label{rmk:5.6}
\emph{
Theorem \ref{thm:1.2} does not exclude the existence of analytic invariant curves
which are not analytic in the perturbation parameter $\e_{0}$. It is an open problem
whether analytic invariant curves may exist under weaker conditions than Bryuno's,
i.e. whether the Bryuno condition is optimal for the system to be analytically conjugated
to a rotation. As observed in Section \ref{sec:1}, Forni's results imply that no such curve
may exist if $\om$ satisfies the condition \eqref{eq:1.5}. However this leaves a gap
of rotation number for which no conclusion can be drawn.
}
\end{rmk}
%%%%%%%%%%%%%%%%%%%%%%%%%%%%%%%%%%%%%%%%%%%%%%%%%

A final comment concerns the Bryuno condition as formulated in Section \ref{sec:1}.
Usually the condition is stated, in any dimension $d$, by requiring for $\oo\in\RRR^{d}$
to satisfy
\begin{equation} \label{eq:5.3}
\gotB_{1}(\oo) = \sum_{n=0}^{\io} \frac{1}{2^{n}} \log \frac{1}{\alpha_{n}(\oo)} < \io , \qquad
\alpha(\oo) = \inf_{\substack{\nn\in\ZZZ^{d} \\ 0<|\nn|\le 2^{n}}} \left| \oo\cdot \nn \right| .
\end{equation}
It is straightforward to check that for $d=1$ and $\oo=(1,\omega)$ the condition
(\ref{eq:6.3}) is equivalent to requiring the Bryuno condition, that is $\gotB_{1}(1,\omega)<\io$
if and only if $\gotB(\om)<\io$ \cite{Y2,Ge3,CM}.

The \emph{Bryuno function} is the solution
of the functional equation \cite{Y1,Y2,MMY1,MMY2}
\begin{equation} \label{eq:5.4}
B(\alpha) = - \log \alpha + \alpha \, B(1/\alpha), \qquad \alpha \in(0,1) .
\end{equation}
Again it is not difficult to show that, if $[\cdot]$ denotes the integer part, one has
$\gotB(\om)<\io$ if and only if $B(\omega-[\omega]) < \io$; see for instance \cite{Y2}.

A somewhat similar Diophantine condition has been proposed by R\"ussmann \cite{R1}
in the context of skew-product systems and later on applied to the study of
KAM tori \cite{P,R2,R3}. Such a condition, known as the \emph{R\"ussmann condition},
can be shown to imply the Bryuno condition \cite{R3}.
However, in dimension $d=2$, is equivalent to the Bryuno condition,
in the sense that, if $\oo=(1,\om)$, then $\omega$ satisfies 
the R\"ussmann condition if and only if $\gotB(\om)<\io$.

%%%%%%%%%%%%%%%%%%%%%%%%%%%%%%%%%%%%%%%%%%%%%%%%%
%%%%%%%%%%%%%%%%%%%%%%%%%%%%%%%%%%%%%%%%%%%%%%%%%
\zerarcounters 
\section{Generalisations}
\label{sec:6} 
%%%%%%%%%%%%%%%%%%%%%%%%%%%%%%%%%%%%%%%%%%%%%%%%%
%%%%%%%%%%%%%%%%%%%%%%%%%%%%%%%%%%%%%%%%%%%%%%%%%

In the general case (\ref{eq:1.1}) the generating function has the form
$S_{\e}(x,x')=S_{0}(x'-x)+\e \s_{\e}(x,x')$, with $S_{0}$ a primitive
of the inverse function of $a$,
which exists under the hypothesis that $a$ satisfies the twist condition.
We discuss explicitly the case in which $\s_{\e}$ does not depend on $\e$ (the
general case can be easily recovered as shown at the end of Section \ref{sec:4}).
Denoting $b=a^{-1}$ and writing the map as
\begin{equation} \label{eq:6.1}
\begin{dcases}
y' = b(x'-x) + \e \, \frac{\partial \sigma_{\e}(x,x')}{\partial x'} , \\
y = b(x'-x) - \e \, \frac{\partial \sigma_{\e}(x,x')}{\partial x} , \\
\end{dcases}
\end{equation}
with $\s_{\e}(x,x')=\ols_{\e}(x,z)=\ols(x,z)$, if we look for a conjugation of the form
\begin{equation} \nonumber
x = \psi + h(\psi) , \qquad y=b(2\pi\om)+H(\psi) ,
\end{equation}
we arrive at the functional equation
\begin{eqnarray} \label{eq:6.2}
& & b \left( 2\pi\om + h(\psi + 2\pi\om) - h(\psi) \right) - b \left( 2\pi\om + h(\psi) - h(\psi-2\pi\om) \right) \nonumber \\
& & \qquad =
\e \, \partial_1 \ols(\psi+h(\psi),2\pi\om+h(\psi+2\pi\om)-h(\psi)) \\
& & \qquad \qquad 
- \, \e \, \partial_2 \ols(\psi+h(\psi),2\pi\om+h(\psi+2\pi\om) -h(\psi)) \nonumber \\
& & \qquad \qquad
+ \, \e \, \partial_2 \ols(\psi-2\pi\om+h(\psi-2\pi\om),2\pi\om+h(\psi)-h(\psi-2\pi\om)) , \nonumber
\end{eqnarray}
which replaces (\ref{eq:2.4}). If we expand
\begin{equation} \label{eq:6.3}
b(2\pi \om+h)=\sum_{k=0}^{\io} B_{k} h^{k} , \qquad B_{k} = \frac{1}{k!} 
\left. \frac{\partial^{k}}{\partial h^{k}} b(2\pi\om + h) \right|_{h=0} ,
\end{equation}
and use that $B_{1} \neq 0$, we obtain
\begin{eqnarray} \label{eq:6.4}
& & h(\psi+2\pi\om) + h(\psi-2\pi\om) - 2 h(\psi) \nonumber \\
& & \quad = \sum_{k=2}^{\io} \frac{B_{k} }{B_1}
\left[ \left( h(\psi) - h(\psi-2\pi\om) \right)^{k} - \left( h(\psi+2\pi\om) - h(\psi) \right)^{k}
\right] \nonumber \\
& & \qquad \qquad + 
\frac{\e}{B_{1}} \, \partial_1 \sigma(\psi+h(\psi),\psi+2\pi\om+h(\psi+2\pi\om)) \\
& & \qquad \qquad \qquad +
\frac{\e}{B_{1}} \, \partial_2 \sigma(\psi-2\pi\om+h(\psi-2\pi\om),\psi+h(\psi)) . \nonumber
\end{eqnarray}
Then by writing $h$ as in (\ref{eq:2.6}), we have, instead of (\ref{eq:2.7}),
\begin{eqnarray} \label{eq:6.5}
& & \hskip-1.truecm \de(\om\nu) \, h^{(k)}_{\nu} =
\sum_{s=2}^{\io} \bar B_{s} 
\sum_{\nu_{1}+\ldots+\nu_{s}=\nu} \sum_{k_{1}+\ldots+k_{s}=k}
\prod_{i=1}^{s} \de_{-}(\om\nu_{i}) h^{(k_{i})}_{\nu_{i}} \nonumber \\
& & \quad -
\sum_{s=2}^{\io} \bar B_{s} 
\sum_{\nu_{1}+\ldots+\nu_{s}=\nu} \sum_{k_{1}+\ldots+k_{s}=k}
\prod_{i=1}^{s} \de_{+}(\om\nu_{i}) h^{(k_{i})}_{\nu_{i}} \nonumber \\
& & \hskip-1.truecm \qquad + 
\sum_{p,q\ge 0} 
\sum_{\substack{\nu_0+\nu_{1}+\ldots+\nu_{p}\\+\mu_{1}+\ldots+\mu_{q}=\nu}}
\sum_{\substack{k_{1}+\ldots+k_{p}\\+k'_{1}+\ldots+k'_{q}=k-1}}
\frac{1}{p!q!} \, (\ii\nu_0)^{p+1} \partial_{z}^{q} \ols_{\nu_0}(2\pi\om) \times \nonumber \\
& & \quad  
\times \Big( \prod_{i=1}^{q} \de_{+}(\om\mu_{i}) \Big)
h^{(k_1)}_{\nu_1}\ldots h^{(k_p)}_{\nu_p}
h^{(k'_1)}_{\mu_1}\ldots h^{(k'_q)}_{\mu_q} \\
& & \hskip-1.truecm \qquad + \sum_{p,q\ge 0} 
\sum_{\substack{\nu_0+\nu_{1}+\ldots+\nu_{p}\\+\mu_{1}+\ldots+\mu_{q}=\nu}}
\sum_{\substack{k_{1}+\ldots+k_{p}\\+k'_{1}+\ldots+k'_{q}=k-1}}
\frac{1}{p!q!} \, (\ii\nu_0)^{p} \partial_{z}^{q+1} \ols_{\nu_0}(2\pi\om) \times \nonumber \\
& & \quad 
\times \left( {\rm e}^{-2\pi\ii \om(\nu_0+\nu_1+\ldots+\nu_p)}
\Big( \prod_{i=1}^{q} \de_{-}(\om\mu_{i}) \Big) -
\Big( \prod_{i=1}^{q} \de_{+}(\om\mu_{i}) \Big) \right)
h^{(k_1)}_{\nu_1}\ldots h^{(k_p)}_{\nu_p}
h^{(k'_1)}_{\mu_1}\ldots h^{(k'_q)}_{\mu_q} , \nonumber
\end{eqnarray}
where we have set $\bar B_{k} = B_{k}/B_{1}$ and $\bar \s_{\nu}(z)=\s_{\nu}(z)/B_{1}$.

By comparing (\ref{eq:6.5}) with (\ref{eq:2.7}), one easily realise 
that a diagrammatic expansion as derived in Section \ref{sec:2} is still possible.
The main difference is that now, to take into account the extra contributions
in the first two lines of (\ref{eq:6.5}), we have to associate with each node
an extra label $\al_{v}=0,1$ with the following properties.
For $\al_{v}=0$ the mode label is $\nu_{v}\in\ZZZ$ and
the corresponding node factor is given by (\ref{eq:2.11}),
with $\sigma_{\nu_{v}}(2\pi\om)$ replaced by $\bar\sigma_{\nu_{v}}(2\pi\om)=
\sigma_{\nu_{v}}(2\pi\om)/B_{1}$. For $\al_{v}=1$ the mode label is $\nu_{v}=0$ and
either $p_{v}=0$ or $q_{v}=0$, while the corresponding node factor is given by
\begin{equation} \label{eq:6.6}
A_{v} = \begin{dcases}
\bar B_{q_{v}} \prod_{\ell \in L_{v}(\theta)^{-}} \de_{-}(\om\nu_{\ell}) , & p_{v} = 0 , \\
- \bar B_{p_{v}}  \prod_{\ell \in L_{v}^{+}(\theta)} \de_{+}(\om\nu_{\ell}) , & q_{v} = 0 .
\end{dcases}
\end{equation}
Then, if we set $V_{\al}(\theta)=\{ v \in V(\theta) : \al_{v}=\al\}$ and
define the order of $\theta$ as $|V_{0}(\theta)|$,
we can define the set $\TT_{k,\nu}$ as the set of trees of order $k$ and
momentum $\nu$ associated with the root line.
The following result is easily proven.

%%%%%%%%%%%%%%%%%%%%%%%%%%%%%%%%%%%%%%%%%%%%%%%%%
\begin{lemma} \label{lem:6.1}
Let $\theta$ be a tree of order $k$. Then the number of nodes of $\theta$ is less
than $2k-1$.
\end{lemma}
%%%%%%%%%%%%%%%%%%%%%%%%%%%%%%%%%%%%%%%%%%%%%%%%%

Apart from these changes the discussion proceeds like in Section \ref{sec:2}.
In particular the analogous of Lemma \ref{lem:2.4} still holds, so that
we have the following result.

%%%%%%%%%%%%%%%%%%%%%%%%%%%%%%%%%%%%%%%%%%%%%%%%%
\begin{lemma} \label{lem:6.2}
One has $\displaystyle{\sum_{\theta\in\TT_{k,0}}\Val(\theta)=0}$ for all $k\ge 1$.
\end{lemma}
%%%%%%%%%%%%%%%%%%%%%%%%%%%%%%%%%%%%%%%%%%%%%%%%%

Of course, the proof requires some extra work, 
because we have to take into account the new nodes: we defer it to Appendix \ref{app:b}.
Then, by reasoning as in Section \ref{sec:2} one proves that the 
perturbation series (\ref{eq:2.6}) is well defined to all orders.

Also the analysis in Sections \ref{sec:3} and \ref{sec:4} can be performed
in the same way, the only difference being that now the derivatives act also
on the factors $\de_{\pm}(\om\nu_{\ell})$ appearing in (\ref{eq:6.6}).
In particular, by defining the set $\gotF(T)$ as in Section \ref{sec:3},
one has that the following analogue of Lemma \ref{lem:3.16} holds true.

%%%%%%%%%%%%%%%%%%%%%%%%%%%%%%%%%%%%%%%%%%%%%%%%%
\begin{lemma} \label{lem:6.3}
For any self-energy cluster $T$ one has
\begin{equation} \nonumber
\sum_{T'\in \gotF(T)} \Val_{T'}(0) = 0 , \qquad
\sum_{T'\in \gotF(T)} \partial_{u} \Val_{T'}(0) = 0 .
\end{equation}
\end{lemma}
%%%%%%%%%%%%%%%%%%%%%%%%%%%%%%%%%%%%%%%%%%%%%%%%%

Again, the proof requires some adaptations with respect to that of Lemma \ref{lem:3.16},
because one has to check that the new node factors do not destroy the
cancellation mechanisms acting on the tree values. We refer once more
to Appendix \ref{app:b} for the details.

Therefore Theorem \ref{thm:1.1} follows. Indeed,
once Lemma \ref{lem:6.3} has been proven, the analysis is essentially a repetition
of the analysis performed in Section \ref{sec:5}, up to minor adaptations, that are required
to take into account the presence of the new factors.
The only relevant novelty with respect to the case (\ref{eq:1.2}) is the discussion,
presented in Appendix \ref{app:b},
of the cancellations leading to Lemma \ref{lem:6.2} and Lemma \ref{lem:6.3}.

%%%%%%%%%%%%%%%%%%%%%%%%%%%%%%%%%%%%%%%%%%%%%%%%%
\vskip.5truecm
\noindent
{\small \textbf{Ackowledgments.}
I'm indebted to Ugo Bessi and Alfonso Sorrentino for useful discussions and
to Livia Corsi for comments and a careful reading of the manuscript.}

%%%%%%%%%%%%%%%%%%%%%%%%%%%%%%%%%%%%%%%%%%%%%%%%%

\appendix

%%%%%%%%%%%%%%%%%%%%%%%%%%%%%%%%%%%%%%%%%%%%%%%%%
%%%%%%%%%%%%%%%%%%%%%%%%%%%%%%%%%%%%%%%%%%%%%%%%%
\zerarcounters 
\section{Cancellations}
\label{app:a} 
%%%%%%%%%%%%%%%%%%%%%%%%%%%%%%%%%%%%%%%%%%%%%%%%%
%%%%%%%%%%%%%%%%%%%%%%%%%%%%%%%%%%%%%%%%%%%%%%%%%

Let $v$ be the node of a tree $\theta$: define $\gotG_{v}$ as the group of
permutations of the subtrees entering $v$, i.e. of the trees that have as the root line
one of the lines entering the node $v$.
Given any line $\ell\in L(\theta)$, define $\gotB_{\ell}$
as the group of transformations which either leave $\beta_{\ell}$ unchanged
or swap it with its opposite;

\vspace{.5truecm}

%%%%%%%%%%%%%%%%%%%%%%%%%%%%%%%%%%%%%%%%%%%%%%%%%
\noindent{\it Proof of Lemma \ref{lem:2.4}}.
When performing the sum over the trees we reason as follows.
Given a tree $\theta$ define $\calF(\theta)$ as the family of trees obtained
through the following operations:
\begin{enumerate}[noitemsep]
\item\label{f1} detach the root line of $\theta$ and re-attach to any other node of $\theta$;
\item\label{f2} apply a transformation in the Cartesian product 
over all $\ell\in L(\theta)$ of the groups $\gotB_{\ell}$;
\item\label{f3} apply a transformation in the Cartesian product
over all $v\in V(\theta)$ of the groups $\gotG_{v}$.
\end{enumerate}

We can rewrite the sum over the trees in $\TT_{k,0}$ as
\begin{equation} \nonumber
\sum_{\theta\in\TT_{k,0}}\Val(\theta)=
\sum_{\theta\in\TT_{k,0}} \frac{1}{|\calF(\theta)|}
\sum_{\theta'\in\calF(\theta)} \Val(\theta') .
\end{equation}
We want to show that for any $\theta\in\TT_{k,0}$ one has
\begin{equation} \label{eq:a.1}
\sum_{\theta'\in\calF(\theta)} \Val(\theta') = 0 .
\end{equation}
For $v'\neq v$ consider the tree $\theta'$ obtained from $\theta$ by detaching the root
line $\ell_{0}$ from the node $v$ and re-attaching it to the node $v'$.
Call $\calP=\calP(v,v')$ the path of lines and nodes connecting $v$ to $v'$:
denote by $w_{1},\ldots,w_{r}$, $r\ge 2$, the ordered nodes in $\calP$, with 
$w_{1}=v$ and $w_{r}=v'$, and by $\ell_{1},\ldots,\ell_{r-1}$ the lines connecting them. 
The node factors and propagators associated with the nodes and lines,
respectively, which do not belong to $\calP$ do not change.

The momenta of the lines $\ell\in\calP$ revert their direction,
so that if the line $\ell$ has momentum $\nu_{\ell}$ in $\theta$,
then it acquires a momentum $\nu_{\ell}'=-\nu_{\ell}$ in $\theta'$
(we are using that the sum over all the mode labels is zero).
In particular the propagators do not change (by parity).
On the contrary, the node factors of the nodes along $\calP$ change
because of the change both of the momenta of the lines and
of the combinatorial factors, as the number of lines entering the nodes may change.

Consider the trees obtained from $\theta$ by applying the transformations in the groups
$\gotB_{\ell}$ for all $\ell\in\calP$, and sum together  all the corresponding values.
We obtain a quantity $A(\theta)$ of the form
\begin{equation} \label{eq:a.2}
A(\theta)= \ii\nu_{v} \, A_{0}(\theta) \, B(\theta) ,
\end{equation}
where 
\begin{eqnarray} \label{eq:a.3}
& & \hskip-1.truecm
A_{0}(\theta) = 
\Biggl[\prod_{j=1}^{r-1} \Biggl( 
\frac{(\ii\nu_{w_{j}})^{P_{j}+1}}{(P_{j}+1)!} \frac{\partial_{z_{w_j}}^{Q_{j}}}{Q_{j}!} +
\left( {\rm e}^{2\pi\ii\om\nu_{\ell_{j}}} - 1 \right) 
\frac{(\ii\nu_{w_{j}})^{P_{j}}}{P_{j}!} \frac{\partial_{z_{w_j}}^{Q_{j}+1}}{(Q_{j}+1)!} \Biggr) \Biggr] \times \\
& & \hskip-1.truecm  \qquad \qquad \qquad \times 
\Biggl[ \frac{(\ii\nu_{w_{r}})^{P_{r}}}{P_{r}!} \frac{\partial_{z_{w_r}}^{Q_{r}}}{Q_{r}!} \Biggr]
\Biggl[ \prod_{j=2}^{r}
\left( \ii\nu_{w_{j}} + \left( {\rm e}^{-2\pi\ii\om\nu_{\ell_{j-1}}} - 1 \right) \partial_{z_{w_j}} \right) \Biggr] 
\prod_{j=1}^{r} \s_{\n_{j}}(2\pi\om) ,
\nonumber
\end{eqnarray}
and $B(\theta)$ takes into account all the other factors.
By writing $\partial_{z_{w_j}}$ we mean that the derivative $\partial_{z}$ acts on $\s_{\nu_{j}}(z)$
and is computed at $z=2\pi\om$.
In (\ref{eq:a.3}) we have called $P_{j}$ and $Q_{j}$ the number of lines
entering the node  $w_{j}$ and not lying on the path $\calP$
which carry a label $+$ and $-$, respectively. This means that,
if we set $p_{j}=p_{w_{j}}$ and $q_{j}=q_{w_{j}}$, one has
$(p_{j},q_{j})=(P_{j}+1,Q_{j})$ if $\be_{\ell_{j}}=+$ and
$(p_{j},q_{j})=(P_{j},Q_{j}+1)$ if $\be_{\ell_{j}}=-$,
while $(p_{r},q_{r})=(P_{r},Q_{r})$. The factor 
$\ii\nu_{v}$ in (\ref{eq:a.2}) is due to the fact that
$\ell_{v}$ is the root line of $\theta$ and hence $\nu_{\ell_{v}}=0$
as $\theta\in\TT_{k,0}$.

When considering the analogous quantity $A(\theta')$ for $\theta'$ one finds 
$A(\theta')=\ii\nu_{v'}A_{0}(\theta')B(\theta')$, 
where we have used that $\nu_{\ell_{v'}}=0$ in $\theta'$,
with $A_{0}(\theta')$ and $B(\theta')$ defined in a similar way,
by taking into account the different orientation of the lines along $\calP$.
However, one has $B(\theta')=B(\theta)$ by construction.

The quantity
$A_{0}(\theta')$ differs from $A_{0}(\theta)$ because the sets $L_{w_{j}}^{\pm}(\theta')$
may have changed with respect to $L_{w_{j}}^{\pm}(\theta)$ as a consequence
of the shifting of the root line from $v$ to $v'$.
Consider the group $\gotG_{w_{1}}$ and identifies the trees which are obtained
from each other by the action of a transformation of the group.
Suppose that $\be_{\ell_{1}}=+$. If $\theta_{1}$ is the subtree with root line $\ell_{1}$
and $s_{1}\ge1$ is the number of subtrees equivalent to $\theta_{1}$ entering $w_{1}$,
the corresponding combinatorial factor is modified into
\begin{equation} \label{eq:a.4}
\frac{1}{p_{1}!} {p_{1} \choose s_{1}} = \frac{1}{s_{1}!(p_{1}-s_{1})!} ,
\end{equation}
while the corresponding factor in $\theta'$ is
\begin{equation} \nonumber
\frac{1}{(p_{1}-1)!} {p_{1}-1\choose s_{1}-1} =
\frac{1}{(s_{1}-1)!(p_{1}-s_{1})!} ,
\end{equation}
which we have to divide by $s_{1}$ to avoid overcounting when considering
all trees in $\cal{F}(\theta)$: indeed, there are in $\theta$ other $s_{1}$ nodes in the
$s_{1}$ trees equivalent to $\theta_{1}$ which produce trees to be identified with $\theta'$.
So that a combinatorial factor (\ref{eq:a.4}) is obtained for the node $w_{1}$ in both cases.
Similarly one deal with the case in which $\be_{\ell_{1}}=-$ and $r_{j}\ge 0$
is the number of subtrees entering $w_{1}$ equivalent to $\theta_{1}$
(of course either $r_{1}$ or $s_{1}$ must be 0), so that we conclude that
the overall combinatorial factor associated to the node $w_{1}$ is
\begin{equation} \label{eq:a.5}
\frac{1}{s_{1}!(p_{1}-s_{1})!} \frac{1}{r_{1}!(q_{1}-r_{1})!}
\end{equation}
for both sets of trees obtained from $\theta$ and $\theta'$.

Analogous considerations hold for the node $w_{r}$, with the roles of $\theta$
and $\theta'$ exchanged. Also the combinatorial factors of the
other nodes along $\cal{P}$ can be dealt with in a similar way: it may happen
that the number of subtrees equivalent to $\theta_{j}$ entering the node $w_{j}$, 
with $1<j<r$, changes either from $s_{j}$ to $s_{j}+1$ or from $r_{j}$ to $r_{j}+1$,
but then one has to divide the corresponding combinatorial factor by
$s_{j}+1$ or $r_{j}+1$, in the respective cases, to avoid overcounting.
On the other hand if the combinatorial factor changes either from $s_{j}$ to $s_{j}-1$ 
or from $r_{j}$ to $r_{j}-1$, then it is the combinatorial factor in $\theta$ 
which must be divided by $s_{j}$ or $r_{j}$, respectively.

By summing together the values of all the inequivalent trees obtained by
following the prescription above, we find that (\ref{eq:a.3}) is replaced by
\begin{eqnarray} \label{eq:a.6}
& & \hskip-1.truecm 
A_{1}(\theta):= \Biggl[ \prod_{j=1}^{r} \Biggl(
\frac{(\ii\nu_{w_{j}})^{P_{j}}}{s_{j}!(P_{j}-s_{j})!} 
\frac{\partial_{z_{w_j}}^{Q_{j}}}{r_{j}!(Q_{j}-r_{j})!} \Biggr) \Biggr] \times \\
& & \hskip-1.truecm  \qquad \times \Biggl[ \prod_{j=1}^{r-1} 
\left( \ii\nu_{w_{j}} + \left( {\rm e}^{2\pi\ii\om\nu_{\ell_{j}}} - 1 \right) \right) \partial_{z_{w_j}} \Biggr]
\Biggl[ \prod_{j=2}^{r}
\left( \ii\nu_{w_{j}} + \left( {\rm e}^{-2\pi\ii\om\nu_{\ell_{j-1}}} - 1 \right) \right) \partial_{z_{w_j}} \Biggr]
\prod_{j=1}^{r} \s_{\n_{j}}(2\pi\om) .
\nonumber
\end{eqnarray}
As noted before, the ordering of the nodes along $\calP$ changes and
each momentum $\nu_{\ell_{j}}$ is replaced with $\nu_{\ell_{j}}'=-\nu_{\ell_{j}}$:
thus, one finds $A_{1}(\theta')=A_{1}(\theta)$ too. In conclusion, if one defines
$\Val_{0}(\theta)=A_{1}(\theta)\,B(\theta)$, one has
\begin{equation} \nonumber
\sum_{\theta'\in\calF(\theta)} \Val(\theta') =
\Val_{0}(\theta) \sum_{v\in V(\theta)} \ii\nu_{v} ,
\end{equation}
which yields (\ref{eq:a.1}).\qed
%%%%%%%%%%%%%%%%%%%%%%%%%%%%%%%%%%%%%%%%%%%%%%%%%

\vspace{.4truecm}

%%%%%%%%%%%%%%%%%%%%%%%%%%%%%%%%%%%%%%%%%%%%%%%%%
\noindent{\it Proof of Lemma \ref{lem:3.16}.}
The argument is very close to that used in the proof of Lemma \ref{lem:2.4}.
Given a self-energy cluster $T$ call $\calF(T)$ the set of all
self-energy clusters $T'$ obtained from $T$ by the following operations:
\begin{enumerate}[noitemsep]
\item\label{sec1} detach the exiting line of $T$ and re-attach it to any $v\in V(\mathring{T})$. 
\item\label{sec2} apply a transformation in the Cartesian product 
over all $\ell\in L(\mathring{T})$ of the groups $\gotB_{\ell}$;
\item\label{sec3} apply a transformation in the Cartesian product
over all $v\in V(\mathring{T})$ of the groups $\gotG_{v}$.
\end{enumerate}

Any line $\ell\in \calP_{T}$ divides $T$ into two disjoint sets $T_{1}(\ell)$ and $T_{2}(\ell)$
such that $T_{2}(\ell)$ contains all the nodes and lines which precede $\ell$ and
$T_{1}(\ell)=T\setminus(\{\ell\}\cup T_{2}(\ell))$. Define $L(T_{j}(\ell))$ and $V(T_{j}(\ell))$,
$j=1,2$, in the obvious way. Call $\calF_{1}(T;\ell)$ the family of
self-energy clusters $T'$ obtained from $T$ by the following operations:
\begin{enumerate}[noitemsep]
\item\label{sec4} attach the exiting line to any node $v\in V(\mathring{T}_{1}(\ell))$ 
and the entering line to any node $w\in V(\mathring{T}_{2}(\ell))$, 
\item\label{sec5} apply a transformation in the Cartesian product 
over all $\ell\in L(\mathring{T})$ of the groups $\gotB_{\ell}$;
\item\label{sec6} apply a transformation in the Cartesian product
over all $v\in V(\mathring{T})$ of the groups $\gotG_{v}$.
\end{enumerate}

Analogously, call $\calF_{2}(T;\ell)$ the family of
self-energy clusters $T'$ obtained from $T$ by the same
operations as $\calF_{1}(T;\ell)$, but by attaching the exiting line 
to any node $v\in V(\mathring{T}_{2}(\ell))$ and the entering line 
to any node $w\in V(\mathring{T}_{1}(\ell))$, instead of applying the operations \ref{sec4}.

We start by proving the first identity of Lemma \ref{lem:3.16}.
Consider the self-energy cluster $T'\in\calF(T)$ obtained from $T$
by detaching the root line from the node $v\in V(T)$ and re-attaching to
the node $v'\in V(\mathring{T})$. Call $\calP=\calP(v,v')$
the path connecting $v'$ to $v$: $\calP$ contains $r$ nodes 
$w_{1}=v \succ w_{2} \succ \ldots \succ w_{r}=v'$, connected by lines $\ell_{1},\ldots,\ell_{r-1}$.

The propagators and the node factors associated with the lines and nodes,
respectively, not belonging to $\calP$ do not change. If we apply
the transformations in the groups $\gotB_{\ell}$ to all lines $\ell\in\calP$ and we sum together
the values of all the self-energy clusters we get,
we realise that the sum of the products of the node factors 
associated with the nodes in $\calP$, is of the form 
$A(T)=\left( \ii\nu_{w_{1}}+({\rm e}^{-2\pi\ii\om\nu} -1)\partial_{z_1}\right) A_{0}(T)\,B(T)$, where
\begin{eqnarray} \nonumber
& & \hskip-1.truecm
A_{0}(T) = \prod_{j=1}^{r-1} \Biggl[
\frac{(\ii\nu_{w_{j}})^{P_{j}+1}}{(P_{j}+1)!} \frac{\partial_{z_{w_j}}^{Q_{j}}}{Q_{j}!} +
\left( {\rm e}^{2\pi\ii\om\nu_{\ell_{j}}} -1 \right)
\frac{(\ii\nu_{w_{j}})^{P_{j}}}{P_{j}!} \frac{\partial_{z_{w_j}}^{Q_{j}+1}}{(Q_{j}+1)!} \Biggr] \times \\
& & \hskip-1.truecm  \qquad \qquad \qquad \times 
\Biggl[ \frac{(\ii\nu_{w_{r}})^{P_{r}}}{P_{r}!} \frac{\partial_{z_{w_r}}^{Q_{r}}}{Q_{r}!} \Biggr]
\prod_{j=2}^{r} \left( \ii\nu_{w_{j}} + \left( {\rm e}^{-2\pi\ii\om\nu_{\ell_{j-1}}} -1 \right) \partial_{z_{w_j}} \right) 
\prod_{j=1}^{r} \s_{\n_{j}}(2\pi\om) , \nonumber
\end{eqnarray}
with $\nu_{\ell_{j}}=\nu_{\ell_{j}}^{0}$ if $\ell_{j}\notin\calP_{T}$ and
$\nu_{\ell_{j}}=\nu_{\ell_{j}}^{0}+\nu_{\ell_{T}'}$ otherwise, and $B(T)$ takes
into account all the other factors. Here we are using the same notations as
described after (\ref{eq:2.3}).
Then we apply the transformations in the groups $\gotG_{w_{j}}$ for $j=1,\ldots,r$,
and reason as in the proof of Lemma \ref{lem:2.4}: summing
together all the inequivalent self-energy clusters produces a quantity of the form
$\left( \ii\nu_{w_{1}}+({\rm e}^{-2\pi\ii\om\nu} -1)\partial_{z_1}\right) A_{1}(T)\,B(T)$, with
\begin{eqnarray} \nonumber
& & \hskip-1.truecm 
A_{1}(T):=\prod_{j=1}^{r} \Biggl[
\frac{(\ii\nu_{w_{j}})^{P_{j}}}{s_{j}!(P_{j}-s_{j})!} 
\frac{\partial_{z_{w_j}}^{Q_{j}}}{r_{j}!(Q_{j}-r_{j})!} \Biggr] \times \\
& & \hskip-1.truecm  \qquad \times \Biggl[ \prod_{j=1}^{r-1} 
\left( \ii\nu_{w_{j}} + \left( {\rm e}^{2\pi\ii\om\nu_{\ell_{j}}} - 1 \right) \partial_{z_{w_j}} \right) \Biggr]
\Biggl[ \prod_{j=2}^{r}
\left( \ii\nu_{w_{j}}+ \left( {\rm e}^{-2\pi\ii\om\nu_{\ell_{j-1}}} - 1 \right) \partial_{z_{w_j}} \right) \Biggr] 
\prod_{j=1}^{r} \s_{\n_{j}}(2\pi\om) , \nonumber
\end{eqnarray}
The same procedure as above, when applied to the self-energy cluster $T'$, produces a factor
$\left( \ii\nu_{w_{r}}+({\rm e}^{-2\pi\ii\om\nu} -1) \partial_{z_{w_r}} \right) A_{1}(T')\,B(T')$,
where $B(T')=B(T)$ and 
\begin{eqnarray} \nonumber
& & \hskip-1.truecm 
A_{1}(T'):=\prod_{j=1}^{r} \Biggl[
\frac{(\ii\nu_{w_{j}})^{P_{j}}}{s_{j}!(P_{j}-s_{j})!} 
\frac{\partial_{z_{w_j}}^{Q_{j}}}{r_{j}!(Q_{j}-r_{j})!} \Biggr] \times \\
& & \hskip-1.truecm  \qquad \times \Biggl[ \prod_{j=2}^{r} 
\left( \ii\nu_{w_{j}}+\left( {\rm e}^{2\pi\ii\om\nu_{\ell_{j-1}}'} - 1 \right) \partial_{z_{w_j}} \right) \Biggr]
\Biggl[ \prod_{j=1}^{r-1}
\left( \ii\nu_{w_{j}}+ \left( {\rm e}^{-2\pi\ii\om \nu_{\ell_{j}}'} -1 \right) \partial_{z_{w_j}} \right) \Biggr]
\prod_{j=1}^{r} \s_{\n_{j}}(2\pi\om) , \nonumber
\end{eqnarray}
where $\nu_{\ell_{j}}'=-\nu_{\ell_{j}}^{0}$ if $\ell_{j}\notin\calP_{T}$ and
$\nu_{\ell_{j}}'=-\nu_{\ell_{j}}^{0}+\nu_{\ell_{T}'}$ otherwise.
Therefore for $\om\nu_{\ell_{T}'}=0$ (so that $\om\nu=0$ as well)
one has $A_{1}(T) = A_{2}(T')$ and hence
\begin{equation} \nonumber
\sum_{T'\in\calF(T)} \Val_{T'}(0)= A_{1}(T)\,B(T) \sum_{v\in V(\mathring{T})} \ii\nu_{v} = 0 ,
\end{equation}
and since all the self-energy clusters in $\calF(T)$ belong to the set $\gotF(T)$
the first identity follows.

Now we pass to the second identity of Lemma \ref{lem:3.16}. First of all,
we note that the derivative $\partial_{u}$ acts on the propagators $\de(\om\nu_{\ell})$
and on the factors ${\rm e}^{\pm2\pi\ii\om\nu_{\ell}}$ of the node factors
corresponding to lines $\ell\in\calP_{T}$.

Given a self-energy cluster $T$ let $\ell$ be a line along the path $\calP_{T}$.
Call $\calA_{T}(u;\ell)$ the contribution
to $\partial_{u}\Val_{T}(u)$ obtained by differentiating the propagator of the line $\ell$.
The differentiated propagator is
\begin{equation} \nonumber
\left. \partial_{u} \calG_{\ell} \right|_{u=0} 
= \left. \partial_{u}\left(1/\de(\om\nu_{\ell}+u) \right) \right|_{u=0} 
= \frac{2\pi \sin 2\pi\om\nu_{\ell}}{\de^{2}(\om\nu_{\ell})} =
\frac{\pi \sin 2\pi\om\nu_{\ell}}{2(\cos 2\pi\om\nu_{\ell} - 1)^{2}}. 
\end{equation}
Consider all the self-energy clusters $T'\in\calF_{1}(T;\ell)$: by reasoning as in the proof
of the first identity one finds that
\begin{equation} \label{eq:3.16}
\sum_{T'\in\calF_{1}(T;\ell)} \calA_{T'}(0;\ell)
= \left. \partial_{u} \calG_{\ell} \right|_{u=0} 
C_{0}(T) \sum_{v\in V(\mathring{T}_{1})} \ii\nu_{v} \sum_{w\in V(\mathring{T}_{2})} \ii\nu_{w}  , 
\end{equation}
for some common factor $C_{0}(T)$. Analogously one finds
\begin{equation} \label{eq:3.17}
\sum_{T'\in\calF_{2}(T;\ell)} \calA_{T'}(0;\ell)
= \left. \partial_{u} \calG_{\ell} \right|_{u=0} C_{0}(T)
\sum_{v\in V(\mathring{T}_{2})} \ii\nu_{v} \sum_{w\in V(\mathring{T}_{1})} \ii\nu_{w} ,
\end{equation}
with the same factor $C_{0}(T)$ as in (\ref{eq:3.16});
we have used that if $\nu=\nu_{\ell_{T}'}$ and $\nu_{\ell'}+\s\nu$,
with $\s\in\{0,1\}$, is the momentum associated with
the line $\ell'\in\calP_{T'}$ in $T'\in\calF_{1}(T;\ell)$, then the momentum
associated with the same line $\ell'$ in $T''\in\calF_{2}(T;\ell)$ is
either $\nu_{\ell'}+\s'\nu$ or $-\nu_{\ell'}+\s''\nu$, with $\s'\,\s''\in\{0,1\}$;
moreover all propagators and node factors have to computed at $u=\om\nu=0$.
In particular the differentiated propagators $\left. \partial_{u} \calG_{\ell} \right|_{u=0}$
in (\ref{eq:3.16}) and (\ref{eq:3.17}) are opposite to each other.
Therefore if we sum together the two contributions
(\ref{eq:3.16}) and (\ref{eq:3.17}) they cancel out.

Call $\calB_{T}(u;\ell)$ and $\calC_{T}(u;\ell)$ the contributions obtained by differentiating the factor 
$({\rm e}^{2\pi\ii\om\nu_{\ell}}-1)\partial_{z_{w(\ell)}}$
and the factor $({\rm e}^{-2\pi\ii\om\nu_{\ell}}-1)\partial_{z_{w'(\ell)}}$, respectively.
The derivative, computed at $\om\nu=0$, gives a factor
$2\pi\ii\,{\rm e}^{2\pi\ii\om\nu_{\ell}}\partial_{z_{w(\ell)}}$ and
$-2\pi\ii\,{\rm e}^{-2\pi\ii\om\nu_{\ell}}\partial_{z_{w'(\ell)}}$, respectively.
Therefore, with respect to $\Val_{T}(0)$, $\calB_{T}(0;\ell)$ and $\calC_{T}(0;\ell)$ contain a factor
\begin{eqnarray}
& & - 2\pi \, {\rm e}^{2\pi\ii\om\nu_{\ell}} \partial_{z_{w(\ell)}} \n_{w'(\ell)} +
2\pi\ii \left (1 - {\rm e}^{2\pi\ii\om\nu_{\ell}} \right) \partial_{z_{w(\ell)}} \partial_{z_{w'(\ell)}}
\nonumber \\
& & \qquad \qquad \quad \hbox{ and } \quad
2\pi \, {\rm e}^{-2\pi\ii\om\nu_{\ell}} \n_{w(\ell)} \partial_{z_{w'(\ell)}} - 
2\pi\ii \left (1 - {\rm e}^{-2\pi\ii\om\nu_{\ell}} \right) \partial_{z_{w(\ell)}} \partial_{z_{w'(\ell)}} , \nonumber
\end{eqnarray}
respectively, instead of
$(\ii\n_{w(\ell)} + ({\rm e}^{2\pi\ii\om\nu_{\ell}}-1)\partial_{w(\ell)})
\, (\ii\n_{w'(\ell)} + ({\rm e}^{-2\pi\ii\om\nu_{\ell}}-1)\partial_{w'(\ell)})$.
Therefore the contributions containing the terms 
$\pm2\pi\ii \left (1 - {\rm e}^{\pm 2\pi\ii\om\nu_{\ell}} \right) \partial_{z_{w(\ell)}} \partial_{z_{w'(\ell)}}$
cancel out by the same mechanism as in the previous case.
When summing the remaining contributions arising from $\calB_{T'}(0;\ell)$
and $\calC_{T'}(0;\ell)$, for $T'\in\calF_{1}(T,\ell)$ we obtain a quantity which is,
up to the sign, the same as the quantity obtained by summing all contributions
$\calB_{T''}(0;\ell)$ and $\calC_{T''}(0;\ell)$ for $T''\in\calF_{2}(T,\ell)$.
The argument is the same as before: again one uses the parity of the
propagator and the fact that $\nu_{\ell}$ becomes $-\nu_{\ell}$ when passing from $T'$ to $T''$.

Finally one has to consider the contributions obtained by differentiating the factors
$(\ii\n_{v}+({\rm e}^{-2\pi\ii\om\nu}-1)\partial_{v})$ and
$(\ii\n_{w}+({\rm e}^{-2\pi\ii\om\nu}-1)\partial_{w})$. Once more, by parity,
the contributions corresponding to $\calF_{1}(T,\ell)$ and $\calF_{2}(T,\ell)$
are opposite to each other and hence cancel out.
\qed
%%%%%%%%%%%%%%%%%%%%%%%%%%%%%%%%%%%%%%%%%%%%%%%%%

%%%%%%%%%%%%%%%%%%%%%%%%%%%%%%%%%%%%%%%%%%%%%%%%%
%%%%%%%%%%%%%%%%%%%%%%%%%%%%%%%%%%%%%%%%%%%%%%%%%
\zerarcounters 
\section{Cancellations in the general case}
\label{app:b} 
%%%%%%%%%%%%%%%%%%%%%%%%%%%%%%%%%%%%%%%%%%%%%%%%%
%%%%%%%%%%%%%%%%%%%%%%%%%%%%%%%%%%%%%%%%%%%%%%%%%

We use the same notations as in Appendix \ref{app:a}. Moreover we need the
following identities.

%%%%%%%%%%%%%%%%%%%%%%%%%%%%%%%%%%%%%%%%%%%%%%%%%
\begin{lemma} \label{lem:b.1}
Let $x_{1},\ldots,x_{n} \in \RRR$ such that $x_{1}+\ldots+x_{n}=0$. Then
\begin{equation} \nonumber
\prod_{i=1}^{n} \left( 1 - {\rm e}^{- x_{i}} \right) =
\prod_{i=1}^{n} \left( {\rm e}^{x_{i}} - 1 \right) .
\end{equation}
\end{lemma}
%%%%%%%%%%%%%%%%%%%%%%%%%%%%%%%%%%%%%%%%%%%%%%%%%

%%%%%%%%%%%%%%%%%%%%%%%%%%%%%%%%%%%%%%%%%%%%%%%%%
\proof
Write ${\rm e}^{x_{i}}-1 = {\rm e}^{x_{i}}(1-{\rm e}^{-x_{i}})$ in each factor of the
second product.
\qed
%%%%%%%%%%%%%%%%%%%%%%%%%%%%%%%%%%%%%%%%%%%%%%%%%

%%%%%%%%%%%%%%%%%%%%%%%%%%%%%%%%%%%%%%%%%%%%%%%%%
\begin{lemma} \label{lem:b.2}
Let $x_{1},\ldots,x_{n} \in \RRR$ such that $x_{1}+\ldots+x_{n}=0$. Then
\begin{equation} \nonumber
\sum_{i=1}^{n} x_{i} \prod_{\substack{j=1 \\ j\neq i}}^{n} \left( 1 - {\rm e}^{- x_{j}} \right) =
\sum_{i=1}^{n} x_{i} \prod_{\substack{j=1 \\ j\neq i}}^{n} \left( {\rm e}^{x_{j}} - 1 \right) .
\end{equation}
\end{lemma}
%%%%%%%%%%%%%%%%%%%%%%%%%%%%%%%%%%%%%%%%%%%%%%%%%

%%%%%%%%%%%%%%%%%%%%%%%%%%%%%%%%%%%%%%%%%%%%%%%%%
\proof We can rewrite
\begin{equation} \nonumber
\sum_{i=1}^{n} x_{i} \prod_{\substack{j=1 \\ j\neq i}}^{n} \left( 1 - {\rm e}^{- x_{j}} \right) =
\sum_{i=1}^{n} x_{i} \left( 1 - {\rm e}^{- x_{i}} \right)^{-1}
\prod_{j=1}^{n} \left( 1 - {\rm e}^{- x_{j}} \right) 
\end{equation}
and, analogously, using that $x_{1}+\ldots+x_{n}=0$,
\begin{equation} \nonumber
\sum_{i=1}^{n} x_{i} \prod_{\substack{j=1 \\ j\neq i}}^{n} \left( {\rm e}^{x_{j}} - 1 \right) =
\sum_{i=1}^{n} x_{i} \left( {\rm e}^{x_{i}} - 1 \right)^{-1} 
\prod_{j=1}^{n} \left( {\rm e}^{x_{j}} - 1 \right) =
\sum_{i=1}^{n} x_{i} \left( {\rm e}^{x_{i}} - 1 \right)^{-1} 
\prod_{j=1}^{n} \left( 1 - {\rm e}^{- x_{j}} \right) ,
\end{equation}
so that the identity follows if we prove that
\begin{equation} \nonumber
\sum_{i=1}^{n} x_{i} \left( 1- {\rm e}^{-x_{i}} \right)^{-1} =
\sum_{i=1}^{n} x_{i} \left( {\rm e}^{x_{i}} - 1 \right)^{-1} 
\end{equation}
for all $x_{1},\ldots,x_{n}$ such that $x_{1}+\ldots+x_{n}=0$.
By writing $\left( 1- {\rm e}^{-x_{i}} \right)^{-1} =
{\rm e}^{x_{i}} \left( {\rm e}^{x_{i}} - 1 \right)^{-1}$, we arrive at
\begin{equation} \nonumber
\sum_{i=1}^{n} x_{i} {\rm e}^{x_{i}} \left( {\rm e}^{x_{i}} - 1 \right)^{-1} =
\sum_{i=1}^{n} x_{i} \left( {\rm e}^{x_{i}} - 1 \right)^{-1} \qquad \Longrightarrow \qquad
\sum_{i=1}^{n} x_{i} \left( {\rm e}^{x_{i}} - 1 \right) \left( {\rm e}^{x_{i}} - 1 \right)^{-1} =0
\end{equation}
which is trivially satisfied if $x_{1}+\ldots+x_{n}=0$.
\qed
%%%%%%%%%%%%%%%%%%%%%%%%%%%%%%%%%%%%%%%%%%%%%%%%%

\vspace{.5truecm}

%%%%%%%%%%%%%%%%%%%%%%%%%%%%%%%%%%%%%%%%%%%%%%%%%
\noindent{\it Proof of Lemma \ref{lem:6.2}}.
Let $\theta$ be a tree in $\TT_{k,0}$. 
If the root line exits a node $v$ with $\al_{v}=1$
the corresponding value vanishes by Lemma \ref{lem:b.1}.
For any node $v \in V(\theta)$ with $\al_{v}=1$ call $\Lambda_{v}$
the set of lines incident with the nodes $v$ (i.e. either exiting or entering $v$)
and, for $\ell\in \Lambda_{v}$, set $\Lambda_{v}(\ell)=\Lambda_{v} \setminus \{\ell\}$ and
\begin{equation} \nonumber
\Delta(v,\ell) := \prod_{\ell' \in \Lambda_{v}(\ell)} \de_{-}(\om\nu_{\ell'}) -
\prod_{\ell' \in \Lambda_{v}(\ell)} \de_{+}(\om\nu_{\ell'}) .
\end{equation}

For $\theta\in\TT_{k,0}$ any $\ell\in L(\theta)$ divides
$\theta$ into two distinct sets $\theta^1_{\ell}$ and $\theta^2_{\ell}$,
where $\theta^1_{\ell}$ is the set of lines and nodes in $\theta$
which precede $\ell$ and $\theta^2_{\ell}$ is the set of lines and nodes
which do not precede $\ell$. If $V(\theta^i_{\ell})$ and $L(\theta^i_{\ell})$
denote, respectively, the set of nodes and the set of lines in $\theta^i_{\ell}$,
one has $V(\theta)=V(\theta^1_{\ell}) \cup V(\theta^2_{\ell})$ and
$L(\theta^1_{\ell}) \cup \{\ell\} \cup L(\theta^2_{\ell})$. 
For $i=1,2$, if $\{v\in V(\theta^{i}_{\ell}) : \al_{v}=1\}=\emptyset$ set
$S_{\ell}(\theta^{i}_{\ell})=1$, otherwise set
\begin{equation} \label{eq:b.1}
S_{\ell}(\theta^{i}_{\ell}) := \prod_{\substack{ v \in V(\theta^{i}_{\ell}) \\ \al_{v}=1}}
\Delta(v,\ell(v)) ,
\end{equation}
where $\ell(v)\in \Lambda_{v}$ is the line incident with $v$ along the
path connecting $v$ with $\ell$. If $\theta'$ is a subset of $\theta^{i}_{\ell}$
we define $S_{\ell}(\theta')$ as in (\ref{eq:b.1}), with the product restricted to
the nodes $v\in V(\theta')$.

Define the family $\calF(\theta)$ as in the proof of Lemma \ref{lem:2.4} 
in Appendix \ref{app:a}, with the further operations of replacing
$p_{v}$ with $q_{v}$ for each node $v \in V(\theta)$.
Define also the family $\calF^{i}_{\ell}(\theta)$, $i=1,2$, as the family
obtained through the same operations considered for $\calF(\theta)$,
but with the operation 1 applied only to the nodes $v\in V(\theta^i_{\ell})$.
We want to prove that
\begin{equation} \label{eq:b.2}
\sum_{\theta'\in\calF^{i}_{\ell}(\theta)} \Val(\theta') =
\nu(\theta^{i}_{\ell}) \, S_{\ell}(\theta^{1}_{\ell}) \, S_{\ell}(\theta^{2}_{\ell}) \, B(\theta) ,
\qquad i=1,2, 
\end{equation}
where we have defined
\begin{equation} \nonumber
\nu(\theta^{i}_{\ell}) := \sum_{v \in V(\theta^{i}_{\ell})} \nu_{v} ,
\end{equation}
and $B(\theta)$ is a suitable common factor independent of $\theta'$.
Note that $\nu(\theta^{1}_{\ell})=-\nu(\theta^{2}_{\ell})=\nu_{\ell}$.

The identity (\ref{eq:b.2}) immediately implies Lemma \ref{lem:6.2}.
Indeed, assuming (\ref{eq:b.2}), one can write, for any line $\ell\in L(\theta)$,
\begin{equation} \nonumber
\sum_{\theta'\in\calF(\theta)} \Val(\theta') =
\sum_{\theta'\in\calF^{1}_{\ell}(\theta)} \Val(\theta') +
\sum_{\theta'\in\calF^{2}_{\ell}(\theta)} \Val(\theta') =
\left( \nu(\theta^{1}_{\ell})+\nu(\theta^{2}_{\ell}) \right)
B(\theta) S_{\ell}(\theta^{1}_{\ell}) S_{\ell}(\theta^{2}_{\ell}) ,
\end{equation}
which yields the assertion since $\nu(\theta^{1}_{\ell})+\nu(\theta^{2}_{\ell})=0$.

Hence, we are left with the proof of (\ref{eq:b.1}).
This is done is by induction on the number of nodes
$v$ with $\al_{v}=1$ contained in $\theta^{i}_{\ell}$.
If $\theta^{i}_{\ell}$ does not contain any node $v$ with $\al_{v}=1$,
then the assertion follows immediately from the same argument
used in proving Lemma \ref{lem:2.4}. Otherwise let $N \ge 1$
be the number of such nodes contained, say, in $\theta^{2}_{\ell}$.
Let $v$ be the node in $V(\theta^2_{\ell})$ such that
$\al_{v}=1$ and $\al_{w}=0$ for all nodes $w \in V(\theta_0)$,
where $\theta_0$ is the set of nodes and lines which precede $\ell$
and do not precede $\ell_{v}$
If $\ell_{1},\ldots,\ell_{s}$, with $s\ge 2$,
is the number of lines entering $v$, call $\theta_{1},\ldots,\theta_{s}$
the subtrees with root lines $\ell_{1},\ldots,\ell_{s}$, respectively.
If we sum the values of all trees in $\calF(\theta)$ obtained by attaching
the root line to a node $v\in V(\theta_0)$ we obtain a contribution
containing a common factor $B_{1}(\theta)$ times
\begin{equation} \nonumber
\nu(\theta_0) \, \Delta(v,\ell_{v}) \prod_{j=1}^{s} S_{\ell}(\theta_{j})
\end{equation}
by Lemma \ref{lem:b.2}. Analogously,
if we sum the values of all trees in $\calF(\theta)$ obtained by attaching
the root line to a node $v\in V(\theta_i)$, we obtain a contribution
containing the same common factor $B_{1}(\theta)$ times
\begin{equation} \nonumber
\nu(\theta_i) \, \Delta(v,\ell_{i}) \prod_{j=1}^{s} S_{\ell}(\theta_{j}) .
\end{equation}
This follows from the the inductive hypothesis,
noting that $S_{\ell}(\theta_{i})=S_{\ell_{i}}(\theta_i)$.
Therefore, by summing together all such contributions, 
we arrive at a common factor times
\begin{equation} \nonumber
\left( \nu(\theta_0) \Delta(v,\ell_{v}) + \sum_{i=1}^{s} 
\nu(\theta_i) \Delta(v,\ell_{i}) \right) \prod_{j=1}^{s} S_{\ell}(\theta_j) .
\end{equation}
Since $\theta\in\TT_{k,0}$, when the root line is attached to a node $v\in\theta_{i}$, 
the momentum associated with the line $\ell_{v}$ is $-\nu_{\ell_{v}}$.
Then, using that $\nu(\theta_0)=\nu_{\ell}-\nu_{\ell_{v}}$ and that
\begin{equation} \nonumber
-\nu_{\ell_{v}} \Delta(v,\ell_{v}) + \sum_{i=1}^{s} 
\nu(\theta_i) \Delta(v,\ell_{i}) = 0
\end{equation}
by Lemma \ref{lem:6.2}, we obtain
\begin{equation} \nonumber
\sum_{\theta'\in\calF^{2}_{\ell}(\theta)} \Val(\theta') =
\nu_{\ell} \, B_{1}(\theta) \, \Delta(v,\ell(v)) \prod_{i=1}^{s} S_{\ell}(\theta_i) ,
\end{equation}
where we have used that $\ell_{v}=\ell(v)$.
The inductive hypothesis also implies that $B_{1}(\theta)$ contains
a factor $\Delta(w,\ell(w))$ for each node $w$ with $\al_{w}=1$
in both $V(\theta^{1}_{\ell})$ and $\cup_{i=1}^{s} V(\theta_i)$,
so that the identity (\ref{eq:b.2}) is proven for $i=2$.
The case $i=1$ can be discussed in the same way.
\qed
%%%%%%%%%%%%%%%%%%%%%%%%%%%%%%%%%%%%%%%%%%%%%%%%%

\vspace{.5truecm}

%%%%%%%%%%%%%%%%%%%%%%%%%%%%%%%%%%%%%%%%%%%%%%%%%
\noindent{\it Proof of Lemma \ref{lem:6.3}}.
One combines the proof of Lemma \ref{lem:3.16} with the proof of Lemma \ref{lem:6.2}.
The proof of the first identity is similar to the proof of Lemma \ref{lem:6.2}, as well
as the first identity of Lemma \ref{lem:3.16} was obtained by following closely the proof
of Lemma \ref{lem:2.4}. To prove the second identity one reasons again as done
for Lemma \ref{lem:3.16}, by using the identity (\ref{eq:b.2}), with $\ell$ being
the line which is differentiated when the derivative acts on the propagator $\calG_{\ell}$.
\qed
%%%%%%%%%%%%%%%%%%%%%%%%%%%%%%%%%%%%%%%%%%%%%%%%%

%%%%%%%%%%%%%%%%%%%%%%%%%%%%%%%%%%%%%%%%%%%%%%%%%

%%%%%%%%%%%%%%%%%%%%%%%%%%%%%%%%%%%%%%%%%%%%%%%%%
% References 
%%%%%%%%%%%%%%%%%%%%%%%%%%%%%%%%%%%%%%%%%%%%%%%%%
%%%%%%%%%%%%%%%%%%%%%%%%%%%%%%%%%%%%%%%%%%%%%%%%%


\begin{thebibliography}{99} 
 
{\small 

\bibitem{BeGi}{
G. Benettin, A. Giorgilli,
\textit{On the Hamiltonian interpolation of near-to-the-identity symplectic mappings
with application to symplectic integration algorithms},
J. Statist. Phys.
\textbf{74} (1994), no. 5-6, 1117-1143. }

\bibitem{BG}{
A. Berretti, G.Gentile,
\textit{Bryuno function and the standard map},
Comm. Math. Phys.
\textbf{220} (2001), no. 3, 623-656. }

\bibitem{Be}{
U. Bessi,
\textit{An analytic counterexample to the KAM theorem},
Ergodic Theory Dynam. Systems 20 (2000), no. 2, 317-333. }

\bibitem{B}{
A.D. Bryuno,
\textit{Analytic form of differential equations. I},
Tr. Mosk. Mat. Obs.
\textbf{25} (1971), 119-262;
Engl. transl. Trans. Moscow Math. Soc. \textbf{25} (1973), 131-288
\textit{Analytic form of differential equations. II},
Tr. Mosk. Mat. Obs.
\textbf{26} (1972), 199-239;
Engl. transl. Trans. Moscow Math. Soc. \textbf{26} (1974), 199-239. }

\bibitem{BC}{
X. Buff, A. Ch\'eritat,
\textit{The Brjuno function continuously estimates the size of quadratic Siegel disks},
Ann. of Math.
\textbf{164} (2006), no. 1, 265-312. }

\bibitem{Ch}{
B.V. Chirikov,
\textit{A universal instability of many-dimensional oscillator systems},
Phys. Rep.
\textbf{52} (1979), no. 5, 264-379. }

\bibitem{CG}{
L. Corsi, G. Gentile,
\textit{Oscillator synchronisation under arbitrary quasi-periodic forcing},
Comm. Math. Phys.
\textbf{316} (2012), no. 2, 489-529. }

\bibitem{CM}{
C. Chavaudret, S. Marmi,
\textit{Reducibility of quasiperiodic cocycles under a Brjuno-R\"ussmann arithmetical condition},
J. Mod. Dyn.
\textbf{6} (2012), no. 1, 59-78. }

\bibitem{D}{
A.M. Davie,
\textit{The critical function for the semistandard map},
Nonlinearity
\textbf{7} (1994), no. 1, 219-229. }

\bibitem{De}{
A. Deprit,
\textit{Canonical transformations depending on a small parameter},
Celestial Mech.
\textbf{1} (1969), 12-30. }

\bibitem{E}{
L.H. Eliasson,
\textit{Absolutely convergent series expansions for quasi periodic motions},
Math. Phys. Electron. J.
\textbf{2} (1996), Paper 4, 33 pp. (electronic). }

\bibitem{F}{
G. Forni,
\textit{Analytic destruction of invariant circles},
Ergodic Theory Dynam. Systems
\textbf{14} (1994), no. 2, 267-298. }

\bibitem{Ga}{
G. Gallavotti,
\textit{Twistless KAM tori},
Comm. Math. Phys.
\textbf{164} (1994), no. 1, 14-156. }

\bibitem{GBG}{
G. Gallavotti, F. Bonetto, G. Gentile,
\textit{Aspects of ergodic, qualitative and statistical theory of motion},
Springer, Berlin, 2004. }

\bibitem{Ge1}{
G. Gentile,
\textit{Brjuno numbers and dynamical systems},
Frontiers in number theory, physics, and geometry. I, 585-599,
Springer, Berlin, 2006. }

\bibitem{Ge2}{
G. Gentile,
\textit{Resummation of perturbation series and reducibility for Bryuno skew-product flows},
J. Stat. Phys.
\textbf{125} (2006), no. 2, 321-361. }

\bibitem{Ge3}{
G. Gentile,
\textit{Degenerate lower-dimensional tori under the Bryuno condition},
Ergodic Theory Dynam. Systems
\textbf{27} (2007), no. 2, 427-457. }

\bibitem{Ge4}{
G. Gentile,
\textit{Quasiperiodic motions in dynamical systems: review of a renormalization group approach},
J. Math. Phys.
\textbf{51} (2010), no. 1, 015207, 34 pp. }

\bibitem{GM}{
G. Gentile, V. Mastropietro,
\textit{Methods for the analysis of the Lindstedt series for KAM tori and renormalizability in classical mechanics.
A review with some applications},
Rev. Math. Phys.
\textbf{8} (1996), no. 3, 393-444. }

\bibitem{GiM}{
A. Giorgilli, S. Marmi,
\textit{Convergence radius in the Poincar\'e-Siegel problem},
Discrete Contin. Dyn. Syst. Ser. S
\textbf{3} (2010), no. 4, 601-621. }

\bibitem{Go}{
Ch. Gol\'e,
\textit{Symplectic twist maps. Global variational techniques},
Advanced Series in Nonlinear Dynamics 18,
World Scientific Publishing Co., Inc., River Edge, NJ, 2001. }

\bibitem{Gr}{
J.M. Greene,
{\it A method for determining a stochastic transition},
J. Math. Phys.
\textbf{20} (1979), no. 6, 1183--1201. }

\bibitem{KK1}{
H. Koch, S. Koci\'c,
\textit{A renormalization group approach to quasiperiodic motion with Brjuno frequencies},
Ergodic Theory Dynam. Systems 30 (2010), no. 4, 1131-1146. }

\bibitem{KK2}{
H. Koch, S. Koci\'c,
\textit{A renormalization approach to lower-dimensional tori with Brjuno frequency vectors},
J. Differential Equations
\textbf{249} (2010), no. 8, 1986-2004. }

\bibitem{KP}{
S. Kuksin, J. P\"oschel,
\textit{On the inclusion of analytic symplectic maps in analytic Hamiltonian flows and its applications},
Seminar on Dynamical Systems (St. Petersburg, 1991), 96-116, Progr. Nonlinear Differential Equations Appl.12,
Birkh\"auser, Basel, 1994. }

\bibitem{LL}{
A.J. Lichtenberg, M.A. Lieberman,
\textit{Regular and chaotic dynamics},
Applied Mathematical Sciences 38, Springer, New York, 1992. }

\bibitem{Mac}{
R.S. MacKay,
\textit{Exact results for an approximate renormalisation scheme and some predictions for the breakup of invariant tori},
Phys. D \textbf{33} (1988), no. 1-3, 240-265. 
\textit{Erratum},
Phys. D \textbf{36} (1989), no. 3, 358. }

\bibitem{MMY1}{
S. Marmi, P. Moussa, J.-Ch.Yoccoz,
\textit{The Brjuno functions and their regularity properties},
Comm. Math. Phys. 186 (1997), no. 2, 265-293. }

\bibitem{MMY2}{
S. Marmi, P. Moussa, J.-Ch.Yoccoz,
\textit{Some properties of real and complex Brjuno functions}, Frontiers in number theory, physics, and geometry. I, 601Ð623,
Springer, Berlin, 2006. }

\bibitem{MS}{
S. Marmi, J. Stark,
\textit{On the standard map critical function},
Nonlinearity
\textbf{5} (1992), no. 3, 743-761. }

\bibitem{M1}{
J. Moser,
\textit{On invariant curves of area-preserving mappings of an annulus},
Nachr. Akad. Wiss. Gšttingen Math.-Phys. Kl. II
\textbf{1962} (1962), 1-20. }

\bibitem{M2}{
J. Moser,
\textit{Monotone twist mappings and the calculus of variations},
Ergodic Theory Dynam. Systems
\textbf{6} (1986), no. 3, 401-413. }

\bibitem{M3}{
J. Moser,
\textit{Minimal solutions of variational problems on a torus},
Ann. Inst. H. Poincar\'e Anal. Non Lin\'eaire
\textbf{3} (1986), no. 3, 229-272. }

\bibitem{N}{
A.I. Neishtadt,
\textit{The separation of motions in systems with rapidly rotating phase},
Prikl. Mat. Mekh.
\textbf{48} (1984), no. 2, 197-204; Engl. transl.
J. Appl. Math. Mech.
\textbf{48} (1985), no. 2, 133-139. }

\bibitem{P}{
J. P\"oschel,
\textit{On elliptic lower-dimensional tori in Hamiltonian systems},
Math. Z.
\textbf{202} (1989), no. 4, 559-608. }

\bibitem{R1}{
H. R\"ussmann,
\textit{On the one-dimensional Schr\"odinger equation with a quasiperiodic potential},
Nonlinear dynamics (Internat. Conf., New York, 1979), pp. 90-107,
Ann. New York Acad. Sci. 357, New York Acad. Sci., New York, 1980. }

\bibitem{R2}{
H. R\"ussmann,
\textit{On the frequencies of quasi-periodic solutions of analytic nearly integrable Hamiltonian systems},
Seminar on Dynamical Systems (St. Petersburg, 1991), 160-183, Progr. Nonlinear Differential Equations Appl. 12,
Birkh\"auser, Basel, 1994. }

\bibitem{R3}{
H. R\"ussmann,
\textit{Invariant tori in non-degenerate nearly integrable Hamiltonian systems},
Regul. Chaotic Dyn.
\textbf{6} (2001), no. 2, 119-204. }

\bibitem{S}{
W.M. Schmidt,
\textit{Diophantine approximations},
Lectures Notes in Mathematics 785,
Springer, Berlin, 1980. }

\bibitem{SM}{
C.L. Siegel, J. Moser,
\textit{Lectures on celestial mechanics},
Springer, New York-Heidelberg, 1971. }

\bibitem{SZ}{
D. Salamon, E. Zehnder,
\textit{KAM theory in configuration space},
Comment. Math. Helv.
\textbf{64} (1989), no. 1, 84-132. }

\bibitem{Y1}{
J.-Ch. Yoccoz,
\textit{Th\'eor\`eme de Siegel, nombres de Bruno et polyn\^omes quadratiques},
Ast\'erisque
\textbf{231} (1995), 3--88. }

\bibitem{Y2}{
J.-Ch. Yoccoz,
\textit{Analytic linearization of circle diffeomorphisms},
Dynamical systems and small divisors (Cetraro, 1998), 125-173,
Lecture Notes in Mathemarics 1784, Springer, Berlin, 2002. }

} 
 
\end{thebibliography}
\end{document}